\documentclass[12pt]{article}
\usepackage{mathrsfs}
\usepackage{amsmath}
\usepackage{amsfonts}
\usepackage{amssymb, amsmath, cite}
\usepackage{color}

\newtheorem{theorem}{Theorem}
\newtheorem{lemma}{Lemma}
\newtheorem{proposition}{Proposition}

\newtheorem{corollary}{Corollary}

\newcommand\be{\begin{equation}}
\newcommand\ee{\end{equation}}
\newcommand\ber{\begin{eqnarray}}
\newcommand\eer{\end{eqnarray}}
\newcommand\berr{\begin{eqnarray*}}
\newcommand\eerr{\end{eqnarray*}}

\newcommand\re{\mathrm{e}}
\newcommand\ri{\mathrm{i}}
\newcommand{\ud}{\mathrm{d}}
\newcommand{\nm}{\nonumber}

\newcommand{\ito}{\int_{\Omega}}
\newcommand{\itr}{\int_{\mathbb{R}^2}}
 \newcommand\wot  {W^{1, 2}(\mathbb{R}^2)}
\newcommand{\vep}{\varepsilon}

\title{ Existence theorems for non-Abelian Chern--Simons--Higgs vortices with  flavor}

\author{Shouxin Chen$^{a}$, \, Xiaosen Han$^{a,b}$,\, Gustavo Lozano$^c$\, and
F.~A.~Schaposnik$^d$
\\~
{\normalsize\it $^a$Institute of Contemporary Mathematics, School of Mathematics, Henan University}\\{\normalsize\it  Kaifeng, Henan 475004, PR China}\\
{\normalsize\it $^b$Taida Institute for Mathematical Sciences, Center for Advanced Study in Theoretical Science}\\{\normalsize\it    National Taiwan University, Taipei,  Taiwan 10617, ROC}\\~
{\normalsize \it  $^{c}$Departamento de F\'\i sica, FCEyN,
 Pabell\'on 1, Ciudad Universitaria}\\
{\normalsize \it  Universidad de Buenos Aires, 1428, Buenos Aires, Argentina}\\
~
{\normalsize \it $^d$\it Departamento de F\'\i sica, Universidad Nacional de La Plata}\\ {\normalsize \it Instituto de F\'\i sica La Plata and CICBA}\\ {\normalsize\it C.C. 67, 1900 La Plata, Argentina}
}

\date{}

\begin{document}
\maketitle
\begin{abstract}
In this paper we establish the existence of vortex solutions for a Chern--Simons--Higgs model with gauge group $SU(N) \times U(1)$ and flavor $SU(N)$, these symmetries ensuring the existence of genuine non-Abelian vortices through a color-flavor locking. Under a suitable ansatz we reduce the
problem to a $2\times 2$ system of nonlinear elliptic equations with exponential terms. We
study this system over the full plane and over a doubly periodic domain, respectively.
For the planar case we use a variational argument to establish the existence result and
derive the decay estimates of the solutions. Over the doubly periodic domain we show
that the system admits at least two gauge-distinct solutions carrying the same physical
energy by using a constrained minimization approach and the mountain-pass theorem.
In both cases we get the quantized vortex magnetic fluxes and electric charges.

\end{abstract}

\section{Introduction}
\setcounter{equation}{0}
Magnetic vortex configurations were   investigated by Abrikosov \cite{Abrikosov} more than fifty years ago in the context of Ginzburg--Landau theory of superconductivity.
Sixteen years later Nielsen  and  Olesen stressed the relevance to high-energy physics of vortex-line solutions of the Abelian Higgs model in the context of  dual string models \cite{NielsenOlesen}.
 Since then the interest on vortices has continued to grow both in condensed-matter and particle physics.

Very early it was observed that  when the Ginzburg--Landau  free-energy parameter ratio  takes the critical value bordering type-I and type-II superconductivity one can
find first-order equations which are equivalent to the more involved second-order Ginzburg--Landau equations  \cite{HA}.  These first-order equations were rediscovered in
the context of high-energy physics in refs.\cite{Bogomolnyi,deVS} with the Ginzburg--Landau parameters identified with the gauge charge and the symmetry
breaking potential coupling constant of the Abelian Higgs model. A rigorous study of the self-dual equations with the
coupling constant ratio at its critical value was presented in \cite{Taubes} where it was proved that the self-dual vortex solutions
  are uniquely determined by a set of $N$ not necessarily
distinct points in the plane corresponding to the zeros of the Higgs field. Every set
of $N$ points determines exactly one such solution.

Interestingly enough,  as first pointed out in \cite{deVS}, the value of the parameter ratio leading to the first-order equations is precisely the one required to extend the Abelian Higgs model to an ${\cal N}=2$  supersymmetric model in $d=3$ space-time dimensions. This supersymmetry connection opened the way to the computation of some data (like the exact particle spectra in certain gauge theories) at strong coupling even when the full theory is not solvable (see \cite{SY} and references therein).

In $(2+1)$ space-time dimensions the usual Maxwell term of the Abelian Higgs model can be replaced by a Chern--Simons (CS) term \cite{DJT,DJT2} leading to the so-called Abelian Chern--Simons--Higgs theory. Also in this case a  specific choice of parameters and of the Higgs potential
(a sixth-order one)    leads to the first-order self-dual equations \cite{Hong,JackiwWeinberg}. The presence of the CS action   drastically changes the vortex solutions which carry in this case both magnetic flux and electric charge. A major interest in the CS theories
  is closely connected to several problems in planar matter physics. In particular, large scale properties of a Quantum
Hall system can be described in terms of a  CS theory  with the Hall conductivity  related to the inverse of the coefficient of the CS action \cite{FZ}. Also,  the central role played by the CS term in bosonization of massive \cite{Fradkin}  and massless \cite{Moreno,AAY} fermions has also been revealed.

The existence of self dual vortex solutions of the Abelian Chern-Simons model with boundary conditions on $\mathbb{R}^2$ has been establish in refs. \cite{Wang,SYcs1}. The model was also  studied in the case of gauge periodic boundary conditions defined on a periodic cell \cite{hoof,waya}  and in this case it has been proved that there is a critical value of the CS coupling parameter above which vortex solutions do not exist \cite{CaffarelliYang,taran96}.
The existence of a critical value of the parameter is related to the area of the periodic domain and has been observed also in { in Abelian and Yang-Mills-Higgs systems in different compact geometries} \cite{Bradlow, GonzalezArroyo, LoMa}. Let us point that  periodic field configuration are  {particularly relevant in the context of condensed matter systems where vortices appear as a lattice array (the Abrikosov lattice)\cite{Abrikosov}}.

Topologically stable vortex solutions in $SU(N)$ gauge theories with both Yang--Mills and Chern--Simons terms were constructed in \cite{dVS1,dVS2} with  the scalar fields breaking the symmetry to  $Z_N$ and gauge fields restricted to the Cartan subalgebra. Under the same assumptions self-dual equations for the pure  non-Abelian Chern--Simons--Higgs  model   were studied in ref. \cite{LMMS}. A proof of the existence   of vortex solutions with the Cartan restriction was given in \cite{Ycs}.  Because in these models the  vortex magnetic fluxes
turn to be in the  Cartan subalgebra direction the solutions can be seen as the result of an Abelian embedding \cite{SY}.

It is known that adding flavor to the Yang--Mills--Higgs theory one can arrange the symmetry breaking so that some global diagonal combination of color and flavor groups survives. This pattern of symmetry breaking is known as color-flavor locking  procedure  (see \cite{SY} and references therein) and leads to genuine  non-Abelian vortices with an orientational moduli related  to the presence of  the  surviving symmetry subgroup. Existence and uniqueness theorems for the solutions of this model were presented in \cite{LY}.

Within the  color-flavor locking  symmetry breaking pattern referred above
  and an appropriate cylindrically symmetric ansatz,  genuine non-Abelian vortex solutions were constructed numerically for the $SU(N)\times U(1)$ Chern-Simons-Higgs theory with $SU(N)$ flavor in \cite{LMMS}.  It is the purpose of this work to present  rigorous existence theorems for this problem which   can be reduced to a  $2\times 2$  system of nonlinear elliptic equations with exponential terms.

  We shall consider  the model defined in two domains: over the full plane and over a doubly periodic domain. Over the plane,  this type of system (with different symmetry groups)) was studied in \cite{Ycs} in a general form by using the Cholesky decomposition of a positive definite matrix to find a variational structure.  However, for our concrete  $2\times 2$ system, we   find  a more  explicitly  variational structure  than that of \cite{Ycs}. Over the doubly periodic domain,
we find a sufficient condition on the coupling parameter such that the system admits at least two different solutions, which are obtained by using a constrained minimization approach    and the mountain-pass theorem.    These two solutions are necessarily gauge inequivalent. Since they both carry the same electric charge and magnetic flux and are of self-dual type, they have the same energy. This phenomenon is in sharp contrast with that in the classical Abelian Higgs model \cite{jata,taubes2}.

The rest of our paper is organized as follows. In section 2 we introduce the non-Abelian Chern--Simons model with gauge group $SU(N)\times U(1)$, reduce it to a $2\times2$ system of nonlinear elliptic partial differential equations with exponential terms, and state our main   results. In section 3 we  establish the existence result for the planar case and derive the decay estimates of the solutions.  Section 4 is devoted to the existence result for the doubly periodic domain case.

\section{The model and main results}
In this section we follow \cite{LMMS} to derive the non-Abelian Chern--Simons--Higgs self-dual equations and state our main results. We study the bosonic sector of the  $\mathcal{N}=2$  SUSY $SU(N)\times U(1)$ Chern--Simons--Higgs action in $2+1$ dimensions
\berr
\mathcal{S}=\!\!\int\!\!\ud^3 x\left\{\frac{\kappa_1}{2}\epsilon^{\mu\nu\rho}F^0_{\mu\nu}A^0_\rho+\frac{\kappa_2}{2}\epsilon^{\mu\nu\rho}\!\!\left(\!F^I_{\mu\nu}
A^I_\rho-\frac13f^{IJK}A_{\mu}^IA_{\nu}^JA_{\rho}^K\!\right)
+(D_\mu\phi^f)^\dagger(D_\mu\phi^f)-V[\phi, \phi^\dagger]\right\}
\eerr
where $\epsilon^{012}=1$, $g^{00}=1$, $f^{IJK}$ are the constructional constants of $SU(N)$. The covariant derivatives and field strengths are defined as
  \berr
   D_\mu \phi^f_a&=&\partial_\mu\phi^f_a+(A_\mu^{SU(N)})_a^b\phi_b^f+(A_\mu^{U(1)})_a^b\phi_b^f,\\
   A_\mu^{SU(N)}&=&A_\mu^I\tau_I, \quad A_\mu^{U(1)}=A_\mu^0\tau_0,\\
   F^0_{\mu\nu}&=&\partial_{[\mu }A_{\nu^0]}, \quad    F^I_{\mu\nu}=\partial_{[\mu }A_{\nu^I]}+f^{IJK}A_\mu^JA_\nu^K,
  \eerr
where $\tau^0$ and $\tau^I$ are the $U(1)$ and $SU(N)$ generators.

The potential $V[\phi, \phi^\dagger]$ is a sixth order polynomial of the form
\ber
 V[\phi, \phi^\dagger]&=&\frac{1}{16\kappa_1^2N^2}\phi_f^\dagger\phi^f\left(\phi_g^\dagger\phi^g-N\xi\right)^2+
 \frac{1}{4\kappa_2^2N^2}\phi_f^\dagger\tau^I\tau^J\phi^f(\phi_g^\dagger\tau^I\phi^g)(\phi_h^\dagger\tau^I\phi^h)\nm\\
 &&-\frac{1}{4\kappa_1\kappa_2N}(\phi_f^\dagger\tau^I\phi^f)^2\left(\phi_g^\dagger\phi^g-N\xi\right),
\label{potencial}\eer
where  $\mu, \nu, \rho=0, 1, 2$ are Lorentz indices, $I, J, K=1, \dots, N^2-1$ are the $SU(N)$  ``color'' group indices, $\tau_I$ are the anti-Hermitian generators of
$SU(N)$. The  scalar complex  mutiplets have both the gauge  index $a, b, c=1,\dots, N$ and flavor index $f, g, h=1, \dots, N$, and can be written as an $N\times N$ matrix.
The choice of the sixth order potential is dictated by the aim of getting first order self-dual(BPS) equations whose static solutions correspond to a lower bound of the energy.
 As it is well known,  the existence of BPS equations is directly related to to the ${\cal N}=2$ supersymmetry of the theory, with the central charge of the supersymmetry
 algebra related to the topological charge of the corresponding BPS solutions \cite{SY}.  It is precisely supersymmetry which requires the particular choice of the potential which,
   in contrast with the usual forth-order one, has two phases. Indeed, up to gauge transformations, the minima of the potential are given by
\berr
\phi^f = 0  & \hspace{1cm}symmetric~phase\\
\phi^f\phi_f^\dagger = \xi{\rm diag}\{1,\ldots,1\} &  \hspace{1 cm}asymmetric~phase
\eerr
In the asymmetric phase, where the  original gauge symmetry is broken,  topologically non-trivial solutions can be found (see the discussion below).
In what follows we set, without loss of generality, $\xi=1$.

The Euler--Lagrange equations of the theory  are
 \berr
  \kappa_1\epsilon_\mu^{\alpha\beta}F^0_{\alpha\beta}&=&J^0_\mu\equiv \phi_f^\dagger\tau^0D_\mu\phi^f-(D_\mu\phi^f)^\dagger\tau^0\phi^f,\\
  \kappa_2\epsilon_\mu^{\alpha\beta}F^I_{\alpha\beta}&=&J^I_\mu\equiv \phi_f^\dagger\tau^ID_\mu\phi^f-(D_\mu\phi^f)^\dagger\tau^I\phi^f,\\
  D_\mu D^\mu\phi^f&=&\frac{\partial V}{\partial \phi_f^\dagger}.
 \eerr
Using the Gauss law,
 \ber
  \kappa_1F^0_{12}=J_0^0=\phi_f^\dagger\tau^0D_0\phi^f-(D_0\phi^f)^\dagger\tau^0\phi^f, \label{g1}\\
  \kappa_2F^I_{12}=J_0^I=\phi_f^\dagger\tau^ID_0\phi^f-(D_0\phi^f)^\dagger\tau^I\phi^f,\label{g2}
 \eer
one finds that
the energy density is given by
 \berr
 \mathcal{H}=(D_0\phi^f)^\dagger(D_0\phi^f)+(D_i\phi^f)^\dagger(D_i\phi^f)+V[\phi, \phi^\dagger].
 \eerr
One can see that the energy
can be written  as a sum of squares
\berr
H= \int\ud^2x  \mathcal{H}= \int\ud^2x\left\{\left[D_0\phi^f-\ri\epsilon \left(\frac{1}{4\kappa_1N}(\phi_g^\dagger\phi^g)\phi^f-\frac{1}{2\kappa_2}\phi^\dagger_g\tau^I\phi^f\right)\right]^\dagger\right.\\
 \times\left.\left[D_0\phi^f-\ri\epsilon \left(\frac{1}{4\kappa_1N}(\phi_g^\dagger\phi^g)\phi^f-\frac{1}{2\kappa_2}\phi^\dagger_g\tau^I\phi^f\right)\right]+(D_{-\epsilon}\phi^f)^\dagger(D_{-\epsilon}\phi^f)+\epsilon\sqrt{2N}F^0_{12}\right\},
\label{H}
\eerr
 where
 \berr
  D_\epsilon=D_1+\ri\epsilon D_2.
 \eerr
Energy minima are then obtained by solving the BPS equations
 \ber
  &&D_0\phi^f-\ri\epsilon \left(\frac{1}{4\kappa_1N}(\phi_g^\dagger\phi^g)\phi^f-\frac{1}{2\kappa_2}\phi^\dagger_g\tau^I\phi^f\right)=0,\label{g3}\\
  &&D_{-\epsilon}\phi^f=0.\label{g4}
 \eer
  Note that because of the presence of the Chern--Simons term,  magnetic vortices are electrically charged and  their magnetic flux ${\cal F}$
\berr
{\cal F}^0 \equiv \int_{\mathbb{R}^2}  F_{12}^0 dx \, ,\;\;\;\; \;\;\;\;
{\cal F}^I \equiv \int_{\mathbb{R}^2}  F_{12}^I dx
\eerr
 and electric charge ${\cal Q}$
\berr
{\cal Q}^0 \equiv \int_{\mathbb{R}^2} J_0^0 dx \, ,\;\;\;\;  \;\;\;\;
{\cal Q}^I \equiv  \int_{\mathbb{R}^2}  J_0^I dx
\eerr
are related. Indeed, one has from  equations \eqref{g1}--\eqref{g2}
\berr
 {\cal Q}^0 = \kappa_1 {\cal F}^0 \, ,\;\;\;\; \;\;\;\;  {\cal Q}^I = \kappa_2 {\cal F}^I.
\eerr

 Since the BPS equations \eqref{g3}--\eqref{g4} are difficult to deal with directly, we make the following ansatz, which coincides with that in   \cite{LMMS} except that cylindrical symmetry is not assumed
  \ber
   \Phi&=&{\rm diag} \{\phi, \dots, \phi, \phi_N\}, \label{i}\\
    A_0^0&=&\frac{1}{\sqrt{2N}}f_0,\quad  A_i^0=-\sqrt{\frac{2}{N}}f_i, \\
    A_0^I&=&0, A^I_i=0, \quad I=1, \dots, N^2-2,\\
      A_0^{N^2-1}&=&\sqrt{\frac{N-1}{2N}}f_0^{N^2-1}, \quad  A_i^{N^2-1}=\sqrt{\frac{2(N-1)}{N}}f_i^{N^2-1},  \label{ff1}
   \eer
where $\phi$ and $\phi_N$ are complex-valued functions,  $f_0, f_i, f_0^{N^2-1}, f_i^{N^2-1}$ are real functions and we have chosen
\berr
\tau^0 =\frac{\ri}{\sqrt{2N}}{\rm diag\{1, \dots, 1\}} \; , \;\;\;\;\;\; \tau^{N^2-1}=\frac{\ri}{\sqrt{2N(N-1)}}{\rm diag\{1, \dots, 1, 1-N\}}.
\eerr

In \eqref{i} we have written the Higgs fields in terms of an  $N\times N$ matrix $\Phi$  with entries $\Phi_{a f}=\phi_a^f$ where $a$ runs over the gauge group indices and $f$ over the flavor indices. The gauge and flavor groups acts on $\Phi$ according to $\Phi \to U \Phi V$ with $U$ an element of the gauge group and $V$ an element of the flavor group. The
 choice of  ansatz for $\Phi$  produces, with appropriate boundary conditions for $\phi$ and $\phi_N$, the spontaneous breaking of both gauge and flavor symmetries with a surviving diagonal global $SU(N)_{C+F}$ in what is known as a color-flavor locking in the vacuum \cite{SY} that will ensure topological stable solutions. Indeed, for the asymmetric phase the first and third terms in the potential \eqref{potencial}  force $\Phi$ to develop a vacuum expectation value while the second one forces it to be
 diagonal. Such vacuum expectation value is preserved only for transformations in which $U = V^{-1}$, which corresponds to perform a global gauge transformation and a related (inverse) global flavor transformation. The relevant
 homotopy group is then $\Pi_1\left(SU(N)\times U(1)/Z_N\right)$, leading to $Z_N$ non-Abelian vortices. Let us finally note that the choice of a different non-trivial $A_i^I$ component together with the corresponding column permutation in $\Phi$  in ansatz \eqref{i}--\eqref{ff1} leads to other $Z_N$ vortex solutions.

 With the above ansatz, the Gauss law \eqref{g1}--\eqref{g2} and the BPS equations \eqref{g3}--\eqref{g4} can be simplified as
  \ber
  \kappa_1F^0_{12}&=&-\frac{1}{N\sqrt{2N}}\left(\!\!\left[[N-1]\!|\phi|^2\!+|\phi_N|^2\right]\!f_0\!\!+\!
  [N-1]\!\left[|\phi|^2\!-\!|\phi_N|^2\right]\!f_0^{N^2-1}\!\!\right)\!, \label{i19}\\
   \kappa_2F^{N^2-1}_{12}&=&-\frac{\sqrt{N-1}}{N\sqrt{2N}}\left(\left[|\phi|^2-|\phi_N|^2\right]f_0+\left[|\phi|^2+[N-1]|\phi_N|^2\right]f_0^{N^2-1}\right),\label{i20}\\
   (\partial_1-\ri\partial_2)\phi&=&\frac{\ri}{N}\left(\left[f_1-f_1^{N^2-1}\right]-\ri\left[f_2-f_2^{N^2-1}\right]\right)\phi,\label{i21}\\
    (\partial_1-\ri\partial_2)\phi_N&=&\frac{\ri}{N}\left(\left[f_1+[N-1]f_1^{N^2-1}\right]-\ri\left[f_2+[N-1]f_2^{N^2-1}\right]\right)\phi_N,\label{i22}\\
    f_0&=&\frac{1}{2\kappa_1}\left([N-1]|\phi|^2+|\phi_N|^2-N\right),\label{i23}\\
     f_0^{N^2-1}&=&\frac{1}{2\kappa_2}\left(|\phi|^2-|\phi_N|^2\right).\label{i24}
  \eer
Without loss of generality we have chosen the upper sign in parameter $\epsilon=\pm 1$ introduced when we wrote the energy density as a sum of squares. The lower sign will just correspond to vortex magnetic fluxes with opposite sign.

From the equations \eqref{i21}--\eqref{i22} of $\phi$ and $\phi_N$, we see that the zeros of them are at most finite and isolated.  We denote  the zero sets of $\phi$ and $\phi_N$
by \[Z_i=\{p_{is}, s=1, \dots, n_i\}, \quad i=1, 2.\]
Let
\berr
 \partial\equiv\frac12(\partial_1-\ri\partial_2),
\eerr
and note
\[\overline{\partial}=\frac12(\partial_1+\ri\partial_2), \quad  \partial\overline{\partial}=\overline{\partial}\partial=\frac14\Delta,\]
by a direct computation, we obtain
\ber
\Delta\ln|\phi|^2=-\sqrt{\frac{2}{N}}F^0_{12}-\sqrt{\frac{2}{N(N-1)}}F^{N^2-1}_{12},\label{F1}\\
\Delta\ln|\phi_N|^2=-\sqrt{\frac{2}{N}}F^0_{12}+\sqrt{\frac{2(N-1)}{N}}F^{N^2-1}_{12}.\label{F2}
\eer

Then, from \eqref{i19}--\eqref{i20}, \eqref{i23}--\eqref{i24} and \eqref{F1}--\eqref{F2}  we have,  away from the zeroes of the Higgs field
\ber
\Delta\ln|\phi|^2= \frac{1}{4N^2}\left\{\frac{1}{\kappa_1^2}\left([N-1]|\phi|^2+|\phi_N|^2\right)\left([N-1]|\phi|^2+|\phi_N|^2-N\right)\right.\nonumber\\
\left.+\frac{N-1}{\kappa_1\kappa_2}\left(|\phi|^2-|\phi_N|^2\right)^2+\frac{1}{\kappa_1\kappa_2}\left([N-1]|\phi|^2+|\phi_N|^2-N\right)\left(|\phi|^2-|\phi_N|^2\right)\right.\nonumber\\
\left.+\frac{1}{\kappa_2^2}\left(|\phi|^2+[N-1]|\phi_N|^2\right)\left(|\phi|^2-|\phi_N|^2\right)\right\},\label{f3}\\
\Delta\ln|\phi_N|^2= \frac{1}{4N^2}\left\{\frac{1}{\kappa_1^2}\left([N-1]|\phi|^2+|\phi_N|^2\right)\left([N-1]|\phi|^2+|\phi_N|^2-N\right)\right.\nonumber\\
\left.+\frac{N-1}{\kappa_1\kappa_2}\left(|\phi|^2-|\phi_N|^2\right)^2-\frac{N-1}{\kappa_1\kappa_2}\left([N-1]|\phi|^2+|\phi_N|^2-N\right)\left(|\phi|^2-|\phi_N|^2\right)\right.\nonumber\\
\left.-\frac{N-1}{\kappa_2^2}\left(|\phi|^2+[N-1]|\phi_N|^2\right)\left(|\phi|^2-|\phi_N|^2\right)\right\}.\label{f4}
\eer

Let
\[u_1=\ln|\phi|^2, \quad u_2=\ln|\phi_N|^2,\quad \lambda\equiv\frac{1}{4\kappa_1^2}, \quad \kappa\equiv\frac{\kappa_1}{\kappa_2}.\]

We can rewrite the above equations \eqref{f3}--\eqref{f4} as
\ber \Delta u_1&=&\lambda\left\{\frac{1}{N^2}\Big([N-1+\kappa]^2\re^{2u_1}-[\kappa-1]\big(N-[N-2][\kappa-1]\big)\re^{u_1+u_2}\right.\nm\\
&&\left.
+[1-\kappa]\big(1+[N-1]\kappa\big)\re^{2u_2}\Big)-\frac1N\Big([N-1+\kappa]\re^{u_1}+[1-\kappa]\re^{u_2}\Big)\right\}\nm\\
&& +4\pi\sum\limits_{s=1}^{n_1}\delta_{p_{1s}},\label{a1}
\eer
\ber
 \Delta u_2&=&\lambda\left\{\frac{1}{N^2}\Big([N-1][1-\kappa][N-1+\kappa]\re^{2u_1}-[N-1][\kappa-1]\big(2+[N-2]\kappa\big)\re^{u_1+u_2}\right.\nm\\
&&\left.
 +\big(1+[N-1]\kappa\big)^2\re^{2u_2}\Big)-\frac1N\Big([N-1][1-\kappa]\re^{u_1}+\big(1+[N-1]\kappa\big)\re^{u_2}\Big)\right\}\nm\\
   &&+4\pi\sum\limits_{s=1}^{n_2}\delta_{p_{2s}}. \label{a2}
\eer

We are interested in the existence of solutions of  \eqref{a1}--\eqref{a2} for two cases. In the first case, we consider the system   \eqref{a1}--\eqref{a2}  over the plane with the topological boundary conditions
 \be
  u_1\to 0,\quad  u_2\to 0,  \quad  |x|\to +\infty. \label{a7}
 \ee
In the second  case we study the equations   over a doubly periodic domain $\Omega$, governing multiple vortices hosted in $\Omega$ such that the field configurations are subject to the 't Hooft boundary condition
\cite{hoof,waya,yang1} under which periodicity is achieved modulo gauge transformations.

Defining the marix   $K$  as
\ber
 K\equiv \frac1N\begin{pmatrix}
 N-1+\kappa&1-\kappa\\
 (N-1)(1-\kappa)&1+(N-1)\kappa
\end{pmatrix} \label{nn1},
\eer
the system  \eqref{a1}--\eqref{a2}  can be rewritten in a compact form as
 \ber
  \Delta u_i=\lambda\left(\sum\limits_{j=1}^2\sum\limits_{k=1}^2\re^{u_j}K_{ij}\re^{u_k}K_{jk}-\sum\limits_{j=1}^2K_{ij}\re^{u_j}\right)+4\pi\sum\limits_{s=1}^{n_i}\delta_{p_{is}}, \, i=1, 2,\label{ca1}
 \eer
 Thus, our model has equations of motion with the same structure as those  studied by Yang  \cite{Ycs} in connection to   the $SU(N)$ Chern Simons model with matter in the {\em adjoint} representation. Indeed, setting $N=2$ and $\kappa=3$ in \eqref{nn1} gives the same equations to those arising in  the $SU(3)$ model studied in \cite{yang1,nota2,Ycs},
 \ber \Delta u_1&=&\lambda\left(4\re^{2u_1}-\re^{u_1+u_2}-2\re^{2u_2}-2\re^{u_1}+\re^{u_2}\right)+4\pi\sum\limits_{s=1}^{n_1}\delta_{p_{1s}}, \label{a5} \\
 \Delta u_2&=&\lambda\left( -2\re^{2u_1}-\re^{u_1+u_2}+4\re^{2u_2}+\re^{u_1}-2\re^{u_2}\right)+4\pi\sum\limits_{s=1}^{n_2}\delta_{p_{2s}}.\label{a6}
\eer

 As already pointed in \cite{Ycs}, an existence theorem for the system \eqref{ca1} in $\mathbb{R}^2$ can be established for more general matrices $K$ not necessarily connected to a specific $SU(N)$ model with adjoint matter. Our equations provide an explicit realization of this idea.  Notice though, that even if the existence of solutions for the system \eqref{ca1} in $\mathbb{R}^2$ can be established as a result of the theorem shown in ref  \cite{Ycs}, our matrix $K$ does not satisfy the hypothesis used in  \cite{Ycs} to derive the decay estimates (the symmetrization of $K$ is not definite positive). It should be stressed that this decay estimates are relevant to make the connection between the topological charge (the magnetic flux) and the number of zeros of the components of the Higgs fields.
We will then present in this paper a new direct variational argument by using an explicit variational structure to get the solutions.  More importantly, we will be able to get the decay estimates of the solutions and quantized fluxes. In addition, it is still an   open problem   existence of solutions   of the problem \eqref{a1}--\eqref{a2} over a  doubly periodic domain. This motivate us to give a complete analysis of the nonlinear elliptic system \eqref{a1}--\eqref{a2} in both cases.

  Our main results  read  as follows.
 \begin{theorem}\label{th1}
  Consider the equations \eqref{a1}--\eqref{a2} over the full plane subject to the topological boundary condition \eqref{a7}.  For any distribution of
  points $p_{i1}, \dots,  p_{in_i}\in \mathbb{R}^2, \,\, i=1, 2$,   $\kappa>0$,  $\lambda>0$, there exists a solution $(u_1, u_2)$ for the equations
  \eqref{a1}--\eqref{a2}  realizing the boundary condition \eqref{a7}. Moreover,  there hold the following decay estimates:  for any small $\vep\in(0, 1)$,   the solution satisfies
   \ber
     |(N-1)u_1+u_2|^2+|u_1-u_2|^2&\le& C(\vep)\re^{-\sigma_0\sqrt{2\lambda}(1-\vep)|x|}, \label{d1}\\
       |\nabla ([N-1]u_1+u_2)|^2+|\nabla (u_1-u_2)|^2&\le& C(\vep)\re^{-\sigma_0\sqrt{2\lambda}(1-\vep)|x|} \label{d2}
   \eer
  as $|x|$ is sufficiently large,   where $C(\vep)$ is a positive constant depending only on $\vep$, $\sigma_0=\min\{1, \kappa\}$.
 \end{theorem}

 \begin{theorem}\label{th2}
   Consider the equations \eqref{a1}--\eqref{a2} over   a doubly periodic domain $\Omega$ in $\mathbb{R}^2$. For any  given   points $p_{i1}, \dots,  p_{in_i} \in \Omega, i=1, 2$, which need not to  be distinct,
    $\lambda>0$, and $\kappa>1$,  we have the following conclusion:

1.  Every solution $(u_1, u_2)$ of \eqref{a1}--\eqref{a2} satisfies
   \be
  \re^{u_1}<1,\quad \re^{u_2}< 1.\label{b02}
    \ee

2. There is a necessary condition
     \be
     \lambda\ge \frac{16\pi([N-1]n_1+n_2)}{N|\Omega|}. \label{a9'}
     \ee
  for the existence of solutions to the  equations \eqref{a1}--\eqref{a2}.

 3. There exist a positive constant $\lambda_0$ such that when $\lambda>\lambda_0$ the  equations
  \eqref{a1}--\eqref{a2} admit at least  two distinct   solutions  over $\Omega$, one of  which   satisfies the behavior
   \be
   \re^{u_1}\to 1, \quad\re^{u_2}\to 1, \quad as \quad \lambda \to +\infty \label{b04}
   \ee
 pointwise a.e. in $\Omega$ and strongly in $L^p(\Omega)$ for any  $p\ge1$.
\end{theorem}

 \begin{theorem} \label{th3}
In both planar and doubly periodic cases, for the solutions  $(u_1, u_2)$ obtained above,   the vortex magnetic fluxes take the
 quantized form
 \begin{eqnarray}
{\cal F}^{U(1)} &\equiv& \int F_{12}^0\ud x = \frac{4\pi}{\sqrt{2N}}\left([N-1]n_1+n_2\right),\label{d6}\\
 {\cal F}^{SU(N)} &\equiv& \int F_{12}^{N^2-1}\ud x = 4\pi\sqrt{\frac{N-1}{2N}} \left(n_1-n_2\right).\label{d6a}
\end{eqnarray}
 \end{theorem}

 %Before turning to the proof of the theorem let us discuss the connection between the zeros of the Higgs fields
%and the (gauge invariant) magnetic flux ${\cal F}^{U(1)}$ associated to the $U(1)$ gauge sector,
%
%Now,
%  from  the equations \eqref{F1}--\eqref{F2} one has
%\ber
%F_{12}^0 &=& -\frac{1}{\sqrt{2N}} \Delta\left((N-1)\ln|\phi|^2 +|\ln\phi_N|^2\right),\\
%F_{12}^{N^2-1} &=& -\sqrt{\frac{N-1}{2N}} \Delta\left(\ln|\phi|^2 - \ln\phi_N|^2\right).
% \eer
% Then by the equations \eqref{f3}--\eqref{f4},
% we have
%  \ber
%   {\cal F}^{U(1)} = \itr F_{12}^0\ud x&=&-\frac{1}{\sqrt{2N}} \itr\Delta \left((N-1)\ln|\phi|^2+|\ln\phi_N|^2\right)\ud x,\label{d5}\\
%     {\cal F}^{SU(N)} = \itr F_{12}^{N^2-1}\ud x&=&-\sqrt{\frac{N-1}{2N}} \itr\Delta \left(\ln|\phi|^2-|\ln\phi_N|^2\right)\ud x.\label{d5a}
%  \eer
%
%Using the  decay estimates  \eqref{d2} and a direct integration, we have
%    \ber
%   {\cal F}^{U(1)}=\\
%      {\cal F}^{SU(N)}=
%  \eer

\section{Planar case}
\setcounter{equation}{0}
In this section we aim to find solutions of  \eqref{a1}--\eqref{a2} under the topological boundary condition \eqref{a7} and establish the decay rate estimates of the solutions, which allows us to get the quantized fluxes stated in Theorem \ref{th3} for the planar case.   We will use a variational argument  as in  \cite{jata} to  prove   Theorem \ref{th1}.

\subsection{Existence of solutions}
Following \cite{jata}, we introduce at this point the background functions
 \[ u_i^0=-\sum\limits_{s=1}^{n_i}\ln(1+\mu|x-p_{is}|^{-2}), \quad \mu>0, \quad i=1, 2,\]
which satisfies
 \ber
  \Delta u_i^0=4\pi\sum\limits_{s=1}^{n_i}\delta_{p_{is}}-h_i, \label{a8}
 \eer
where
  \be
   h_i=\sum\limits_{s=1}^{n_i}\frac{4\mu}{(\mu+|x-p_{is}|^2)^2}, \quad i=1,2. \label{a9}
  \ee

Writing $u_i=u_i^0+v_i$,   we then  recast \eqref{a1}--\eqref{a2} as
 \ber \Delta v_1\!\!\!&=&\!\!\lambda\left\{\frac{1}{N^2}\Big([N-1+\kappa]^2\re^{2u_1^0+2v_1}-[\kappa-1]\big(N-[N-2][\kappa-1]\big)\re^{u_1^0+u_2^0+v_1+v_2}\right.\nm\\
\!\!&&\!\!\!\left.+[1-\kappa]\big(1+[N-1]\kappa\big)\re^{2u_2^0+2v_2}\Big)-\frac1N\Big([N-1+\kappa]\re^{u_1^0+v_1}+[1-\kappa]\re^{u_2^0+v_2}\Big)\right\}\nm\\
\!\!\!&&\!\!+h_1,\label{a10}\\
 \Delta v_2\!\!\!&=&\!\!\lambda\left\{\frac{[N-1][1-\kappa]}{N^2}\Big([N-1+\kappa]\re^{2u_1^0+2v_1}+\big(2+[N-2]\kappa\big)\re^{u_1^0+u_2^0+v_1+v_2}\right.\nm\\
\!\!\!&&\!\!\left.
\!\! +\frac{\big(1+[N-1]\kappa\big)^2}{N^2} \re^{2u_2^0+2v_2}\!\Big)\!-\!\frac1N\!\Big(\![N-1][1-\kappa]\re^{u_1^0+v_1}\!+
\!\big(1\!+\![N-1]\kappa\big)\re^{u_2^0+v_2}\!\Big)\!\right\}\nm\\
   \!\!\!&&\!\!\!+h_2.\label{a11}
\eer

The boundary condition \eqref{a7} now reads as
 \be
  v_1\to0, \quad v_2\to0, \quad |x|\to +\infty.\label{a7'}
 \ee

To  see the  variational structure \eqref{a10}--\eqref{a11} clearly,
it is convenient to rewrite the equations \eqref{a10}--\eqref{a11} equivalently as
\ber
 &&\left(N-1+\frac1\kappa\right)\Delta v_1+\left(1-\frac1\kappa\right)\Delta v_2=\lambda\left([N-1+\kappa]\re^{2u_1^0+2v_1}\right.\nm\\
 &&\quad \left.-N\re^{u_1^0+v_1}-[\kappa-1]\re^{u_1^0+u_2^0+v_1+v_2}\right)+ \left(N-1+\frac1\kappa\right)h_1+\left(1-\frac1\kappa\right)h_2,\label{a12}\\
  &&\left(1-\frac1\kappa\right)\Delta v_1+\left(\frac{1}{N-1}+\frac1\kappa\right)\Delta v_2=\lambda\left(\left[\frac{1}{N-1}+\kappa\right]\re^{2u_2^0+2v_2}\right.\nm\\
 &&\quad \left.-\frac{N}{N-1}\re^{u_2^0+v_2}-[\kappa-1]\re^{u_1^0+u_2^0+v_1+v_2}\right)+ \left(1-\frac1\kappa\right)h_1+\left(\frac{1}{N-1}+\frac1\kappa\right)h_2.\label{a13}
\eer

We will work on the   space $\wot\times\wot$.
Let us define
\ber
A(N,\kappa)  \equiv \frac{1}{N}\begin{pmatrix}
 N-1+\frac{1}{\kappa} &1-\frac{1}{\kappa}\\
 1-\frac{1}{\kappa}&\frac{1}{N-1}+\frac{1} {\kappa}
\end{pmatrix} \label{11}.
\eer
and introduce $\mathbf{v}^t=(v_1,v_2)$, $\mathbf{q}^t=(e^{u_1^0+v_1}-1,e^{u_2^0+v_2}-1)$, $\mathbf{h}^t=(h_1,h_2)$.
Then, it is straightforward to check that the equations \eqref{a12}--\eqref{a13} are the Euler-Lagrange equations of the following functional
 \be
  I(v_1, v_2)=\itr \frac{1}{2} {\nabla \mathbf{v}^t} A(N,\kappa) {\nabla \mathbf{v}} + \frac{\lambda}{2}
 \mathbf{q}^t A(N,\frac{1}{\kappa}) \mathbf{q}+ \mathbf{h}^t A(N,\kappa) \mathbf{v}
 .\label{a14}
\ee

%It is straightforward to check that the equations \eqref{a12}-\eqref{a13} are the Euler-Lagrange equations of the following functional
% \ber
%  I(v_1, v_2)&=&\frac12\left(N-1+\frac1\kappa\right)\|\nabla v_1\|_2^2+\frac12\left(\frac{1}{N-1}+\frac1\kappa\right)\|\nabla v_2\|_2^2+\left(1-\frac1\kappa\right)\itr\nabla v_1\cdot\nabla %v_2\ud x\nm\\
%   &&+\frac\lambda2\left([N-1+\kappa]\itr[\re^{u_1^0+v_1}-1]^2\ud x+\left[\frac{1}{N-1}+\kappa\right]\itr[\re^{u_2^0+v_2}-1]^2\ud x\right.\nm\\
%   &&\left.+2[1-\kappa]\itr[\re^{u_1^0+v_1}-1][\re^{u_2^0+v_2}-1]\ud x\right)\nm\\
%   &&+\itr\left(\left[N-1+\frac1\kappa\right]h_1+\left[1-\frac1\kappa\right]h_2\right)v_1\ud x\nm\\
%   &&+\itr\left(\left[1-\frac1\kappa\right]h_1+\left[\frac{1}{N-1}+\frac1\kappa\right]h_2\right)v_2\ud x.\label{a14}
%\eer
Then, to solve the equations \eqref{a12}--\eqref{a13} (or equivalently \eqref{a10}--\eqref{a11}), we just need to find the critical points of the functional $I$ defined above.

To seek the critical points of the functional $I$, we first show that it is coercive over $\wot\times\wot$.

 It is easy to see that   $A(N,\kappa)$, as defined in \eqref{11}, is positive definite,   and the smaller eigenvalue is
 \be
 \alpha_0 (\kappa)\equiv\frac12\left(\frac{(N-1)^2+1}{N-1}+\frac2\kappa-\sqrt{\frac{N^2(N-2)^2}{(N-1)^2}+4\left(1-\frac1\kappa\right)^2}\right)>0\label{a140}
 \ee
 for any $N\ge2, \, \kappa>0$.

 Then
 \ber
 I(v_1,v_2) &&\ge\frac{\alpha_0(\kappa)}{2}\left(\|\nabla \mathbf{v}\|_2^2\right) + \alpha_0(\kappa^{-1}) \frac{\lambda}{2}\left(\|\mathbf{q}\|_2^2\right)\label{a14a} +
 \itr \mathbf{h}^t A(M,\kappa) \mathbf{v}
 \eer

 From  the expression of $h_i$ \eqref{a9}, we see that
 \be
  \|h_i\|_2\le \frac{C}{\sqrt{\mu}}, \quad i=1, 2,  \label{a18}
 \ee
 here and in the following  we use $C$ to denote a generic positive constant independent of $\mu$.
Then it follows from the H\"{o}lder inequality and \eqref{a18} that
 \ber
   \itr\left(\left[N-1+\frac1\kappa\right]h_1+\left[1-\frac1\kappa\right]h_2\right)v_1\ud x&\ge& -\frac{C}{\sqrt{\mu}}\|v_1\|_2,\label{a19} \\
   \itr\left(\left[1-\frac1\kappa\right]h_1+\left[\frac{1}{N-1}+\frac1\kappa\right]h_2\right)v_2\ud x &\ge& -\frac{C}{\sqrt{\mu}}\|v_2\|_2.\label{a20}
  \eer

Next we need to control $L^2$-norms in eq. \eqref{a14a}.
  Noting the fact  $\re^{u_i^0}-1\in L^2(\mathbb{R}^2)$  and using the elementary inequality $|\re^t-1|\ge \frac{|t|}{1+|t|}, \quad t\in\mathbb{R}$, we have
  \ber
   \itr\left(\re^{u_i^0+v_i}-1\right)^2\ud x&=&\itr\left(\re^{u_i^0}[\re^{v_i}-1]+\re^{u_i^0}-1\right)^2\ud x\nm\\
    &\ge&\frac12\itr\re^{2u_i^0}\left(\re^{v_i}-1\right)^2 \ud x-\itr\left(\re^{u_i^0}-1\right)^2\ud x\nm\\
     &\ge& \frac12\itr\re^{2u_i^0}\frac{|v_i|^2}{(1+|v_i|)^2}\ud x-C, \quad i=1,2. \label{a21}
  \eer

From the definition of $u_i^0$, we see that $\re^{2u_i^0}$ satisfies  $0\le \re^{2u_i^0}<1$  and  vanishes at the points $p_{i1}, \dots, p_{in_i}, \,\, i=1, 2$.
To proceed further, we use a decomposition of $\mathbb{R}^2$ as in \cite{jata}
 \be
  \mathbb{R}^2=\Omega_1^i\cup\Omega_2^i, \quad i=1, 2,  \label{a22}
 \ee
 where
 \ber
 \Omega_1^i=\left\{x\in \mathbb{R}^2\, \big| \, \re^{2u_i^0}\le \frac12\right\}, \quad   \Omega_2^i=\left\{x\in \mathbb{R}^2\, \big| \, \re^{2u_i^0}\ge \frac12\right\}, \quad i=1, 2.\label{a23}
 \eer

Next we need the inverse H\"{o}lder inequality (see \cite{Wang}).
 \begin{lemma}
  For any   measurable functions $g_1, g_2$ over $\Omega$, there holds
    \be
     \ito|g_1g_2|\ud x \ge \left(\ito|g_1|^q\ud x\right)^{\frac1q} \left(\ito|g_2|^{q'}\ud x\right)^{\frac{1}{q'}}, \label{a24}
    \ee
   where $q, q'\in \mathbb{R}$,  $0<q<1$, $q'<0$ with $\frac1q+\frac{1}{q'}=1$.
 \end{lemma}

On $\Omega_1^i$, we have  $0\le \re^{2u_i^0}\le \frac12$ and  $\re^{2u_i^0}$ approaches $0$ at most $4n_i$ order near the vortex points $p_{is}$, $s=1, \dots, n_i, \, i=1, 2$.
We choose  $q_i'$  to satisfy  $-\frac{1}{2n_i}<q_i'<0$, then the integrals
  \[\int_{\Omega_1^i} \re^{2q_i'u_i^0}\ud x, \quad i=1,2\]
  exist and are positive constants.    In view of the inverse Holder inequality \eqref{a24}, we obtain
  \ber
    \int_{\Omega_1^i}\re^{2u_i^0}\frac{|v_i|^2}{(1+|v_i|)^2}\ud x&\ge& \left(\int_{\Omega_1^i}\frac{|v_i|^{2q_i}}{(1+|v_i|)^{2q_i}}\ud x\right)^{\frac{1}{q_i}}\left(\int_{\Omega_1^i}\re^{2q_i'u_i^0}\ud x\right)^{\frac{1}{q_i'}}\nm\\
     &\ge& C \left(\int_{\Omega_1^i}\frac{|v_i|^{2q_i}}{(1+|v_i|)^{2q_i}}\ud x\right)^{\frac{1}{q_i}}, \quad i=1, 2,\label{a25}
  \eer
  where $0<q_i<\frac{1}{1+2n_i}, \,i=1, 2.$

  Since  $0<\frac{|v_i|}{1+|v_i|}<1, $  using Young's  inequality, we have
 \ber
  \left(\int_{\Omega_1^i}\frac{|v_i|^{2q_i}}{(1+|v_i|)^{2q_i}}\ud x\right)^{\frac{1}{q_i}}&\ge&\left(\int_{\Omega_1^i}\frac{|v_i|^2}{(1+|v_i|)^2}\ud x\right)^{\frac{1}{q_i}}\nm\\
   &\ge&\frac12\int_{\Omega_1^i}\frac{|v_i|^2}{(1+|v_i|)^2}\ud x-C.\label{a26}
  \eer

  Inserting \eqref{a26} into \eqref{a25}, we find that
    \ber
    \int_{\Omega_1^i}\re^{2u_i^0}\frac{|v_i|^2}{(1+|v_i|)^2}\ud x&\ge& C\int_{\Omega_1^i}\frac{|v_i|^2}{(1+|v_i|)^2}\ud x-C, \quad i=1, 2.\label{a27}
   \eer

  On $\Omega_2^i$, we have
     \ber
    \int_{\Omega_2^i}\re^{2u_i^0}\frac{|v_i|^2}{(1+|v_i|)^2}\ud x&\ge& \frac12\int_{\Omega_2^i}\frac{|v_i|^2}{(1+|v_i|)^2}\ud x, \quad i=1, 2.\label{a28}
   \eer

  Combining \eqref{a21}, \eqref{a27} and \eqref{a28}, we finally obtain
    \be
   \|q_i\|_2= \itr\left(\re^{u_i^0+v_i}-1\right)^2\ud x\ge C\itr\frac{|v_i|^2}{(1+|v_i|)^2}\ud x-C, \quad i=1, 2.\label{a29}
    \ee

Let us now analize $\|\nabla v_i\|_2^2$.  Noting the interpolation inequality over $\wot$
   \be
    \itr w^4\ud x\le 2\itr w^2\ud x\itr |\nabla w|^2\ud x, \quad \forall\, w\in \wot, \label{a30}
   \ee
  we obtain
    \ber
   &&\left(\itr|v_i|^2\ud x\right)^2\nm\\
   &&=\left(\itr\frac{|v_i|}{1+|v_i|}[1+|v_i|]|v_i|\ud  x\right)^2\nm\\
    &&\le\itr\frac{|v_i|^2}{(1+|v_i|)^2}\ud x\itr\big(|v_i|+|v_i|^2\big)^2\ud x\nm\\
    &&\le4\itr\frac{|v_i|^2}{(1+|v_i|)^2}\ud x\itr|v_i|^2\ud x\left(\itr|\nabla v_i|^2+1\right)\nm\\
    &&\le\frac12 \left(\itr|v_i|^2\ud x\right)^2+C\left(\left[\itr\frac{|v_i|^2}{[1+|v_i|]^2}\ud x\right]^4+\left[\itr|\nabla v_i|^2\ud x\right]^4+1\right), \label{a31}
  \eer
  which  implies
   \ber
   \|v_i\|_2&\le& C\left(\itr \frac{|v_i|^2}{[1+|v_i|]^2}\ud x+\itr|\nabla v_i|^2\ud x+1\right), \quad i=1, 2. \label{a32}
   \eer

  Then from \eqref{a14a}, \eqref{a19}, \eqref{a20}, and \eqref{a29},   we conclude that
   \ber
    I(v_1, v_2)&\ge&\frac{\alpha_0(\kappa)}{2}\|\nabla \mathbf{v}\|_2^2
      +C\left(\itr \frac{|v_1|^2}{[1+|v_1|]^2}\ud x+\itr \frac{|v_2|^2}{[1+|v_2|]^2}\ud x\right)\nm\\
      &&-\frac{C}{\sqrt{\mu}}\big(\|v_1\|_2+\|v_2\|_2\big)-C.\label{a33}
   \eer

    %\ber
    %I(v_1, v_2)&\ge&\frac{\alpha_0}{2}\left(\|\nabla v_1\|_2^2+\|\nabla v_2\|_2^2\right)
    %  +C\left(\itr \frac{|v_1|^2}{[1+|v_1|]^2}\ud x+\itr \frac{|v_2|^2}{[1+|v_2|]^2}\ud x\right)\nm\\
    %  &&-\frac{C}{\sqrt{\mu}}\big(\|v_1\|_2+\|v_2\|_2\big)-C.\label{a33}
   %\eer

  At this point, taking $\mu$ sufficiently large in \eqref{a33} and using \eqref{a32}, we have
     \ber
    I(v_1, v_2)&\ge&C\left( \|\nabla \mathbf{v}\|_2^2 +\itr \frac{|v_1|^2}{[1+|v_1|]^2}\ud x+\itr \frac{|v_2|^2}{[1+|v_2|]^2}\ud x\right)-C.\label{a34}
   \eer
Thus, from \eqref{a34} and \eqref{a32}, we get
     \ber
    I(v_1, v_2)&\ge&C\left(\|v_1\|_{\wot}+\|v_2\|_{\wot}\right)-C,\label{a35}
   \eer
  which gives the coerciveness of the functional $I$ over $\wot\times \wot$.

 It is easy to see that the functional $I$ is continuous, differentiable, and weakly lower semi-continuous on $\wot\times\wot$. Then by \eqref{a35},
 we infer that the functional $I$ admits a critical point $(v_1, v_2)\in \wot\times\wot$, which is a weak solution of the equations  \eqref{a12}-\eqref{a13}.

Using the well-known inequality
 \berr
 \|\re^{w}-1\|_2\le C\exp(C\|w\|_{\wot}^2), \quad \forall \, w\in \wot,
 \eerr
we see that the  right hand sides of the equations \eqref{a12}-\eqref{a13} belong to $L^2(\mathbb{R}^2)$.  Hence using the $L^2$ elliptic estimates
we conclude that    $(v_1, v_2)\in W^{2,2}(\mathbb{R}^2)\times W^{2,2}(\mathbb{R}^2)$,  which  gives  the desired boundary condition \eqref{a7'} at infinity.
Similarly, we  obtain that  $\partial_iv_1, \partial_iv_2\to 0$ $(i=1, 2)$ when $|x|\to\infty$.

\subsection{ Decay estimates and quantized fluxes}

In this subsection our purpose is   to
derive the decay rates   of $(N-1)u_1+u_2, u_1-u_2$ and their derivatives   when $|x|\to \infty$. As an application of the decay estimates we  may calculate  the quantized fluxes stated  in Theorem \ref{th3} for the planar case.

To establish  the decay estimate, it is convenient to write the equations in a vector form.
We will use the following notation
\[\mathbf{u}=(u_1, u_2)^\tau,\quad \mathrm{U}={\rm diag}\{\re^{u_1}, \re^{u_2}\},\quad \mathbf{U}=(\re^{u_1}, \re^{u_2})^\tau, \quad \mathbf{1}=(1, 1)^\tau.\]
Then away form the vortex points the equations  \eqref{a1}--\eqref{a2} can be rewritten as
 \ber
  \Delta \mathbf{u}=\lambda K\mathrm{U}K(\mathbf{U}-\mathbf{1}),\label{a36}
 \eer
where $K$ is defined by \eqref{nn1}.

Let
\ber
 O\equiv  \begin{pmatrix}
 N-1 &1 \\
 1&-1
\end{pmatrix} \label{n1},\quad
\mathbf{w}=(w_1, w_2)^\tau\equiv O\mathbf{u}.
\eer

 When $|x|>R$ with $R>0$ sufficiently large such that $R>|p_{is}|$ for $s=1, \dots, n_i$ and $i=1, 2$,  we have
  \ber
  \Delta \mathbf{u}=\lambda K\mathrm{U}K\mathrm{U_\xi}\mathbf{u}=\lambda K^2\mathbf{u}+\lambda \big(K\mathrm{U}K\mathrm{U_\xi}\mathbf{u}-K^2\mathbf{u}\big),\label{a37}
 \eer
 where $\mathrm{U_\xi}\equiv{\rm diag}\{\re^{u_1^\xi}, \re^{u_2^\xi}\}$, and  $u_i^\xi$ lies between  $0$ and $u_i$ for $i=1, 2$.
From \eqref{a37}, we see that when $|x|>R$
 \ber
  \Delta \mathbf{w}&=&\lambda OK^2O^{-1}\mathbf{w}+\lambda \big(OK\mathrm{U}K\mathrm{U_\xi}O^{-1}\mathbf{w}-OK^2O^{-1}\mathbf{w}\big)\nm\\
   &=& \lambda D\mathbf{w}+\lambda \big(OK\mathrm{U}K\mathrm{U_\xi}O^{-1}\mathbf{w}-D\mathbf{w}\big),\label{n2}
 \eer
where $D\equiv{\rm diag}\{1, \kappa^2\}$.

 Then as $|x|>R$ we have
 \ber
  \Delta |\mathbf{w}|^2\ge 2\mathbf{w}^\tau\Delta \mathbf{w}&=&2\lambda\mathbf{w}^\tau D\mathbf{w}+\lambda \mathbf{w}^\tau\big(OK\mathrm{U}K\mathrm{U_\xi}O^{-1}\mathbf{w}-D\mathbf{w}\big) \nm\\
  &\ge&2\lambda\sigma_0^2|\mathbf{w}|^2-f(x)|\mathbf{w}|^2, \label{a38}
 \eer
where  $\sigma_0=\min\{1, \kappa \}$,   $f(x)\to 0$ as $|x|\to\infty$.  Therefore, for any sufficiently small $\vep\in(0, 1)$, there exists an $R_\vep>R$ such that
 \ber
   \Delta |\mathbf{w}|^2\ge 2\lambda\sigma_0^2\left(1-\frac\vep2\right)|\mathbf{w}|^2, \quad |x|>R_\vep. \label{a40}
 \eer
Noting that $|\mathbf{w}|^2=0$ at infinity,  we conclude from \eqref{a40} that there exits a positive constant $C(\vep)$ such that
  \ber
   |\mathbf{w}|^2\le C(\vep)\re^{-\sigma_0\sqrt{2\lambda}(1-\vep)|x|}, \quad \quad |x|>R_\vep, \label{a41}
  \eer
which gives the desired estimate \eqref{d1}.

To get the decay estimate for the derivatives of $w_1$ and $w_2$, we follow the same procedure.  Let $\partial$ be any one of the two partial derivatives $\partial_1$ and $\partial_2$.
Then from \eqref{a36} we get
 \ber
  \Delta (\partial\mathbf{u})&=&\lambda K\mathrm{U}\mathrm{V}K(\mathbf{U}-\mathrm{1})+\lambda K\mathrm{U}K\mathrm{U}\partial\mathbf{u}\nm\\
  &=& \lambda K^2\partial\mathbf{u} + \lambda\big( K\mathrm{U}\mathrm{V}K(\mathbf{U}-\mathrm{1})+K\mathrm{U}K\mathrm{U}\partial\mathbf{u}-K^2\partial\mathbf{u}\big),\label{a42}
 \eer
where $ \mathrm{V}\equiv{\rm diag}\{\partial u_1, \partial u_2\}$.

Hence we have
 \ber
  \Delta (\partial\mathbf{w})&=& \lambda OK^2O^{-1}\partial\mathbf{w}   + \lambda\big( OK\mathrm{U}\mathrm{V}K(\mathbf{U}-\mathrm{1})+OK\mathrm{U}K\mathrm{U}O^{-1}\partial\mathbf{w}-OK^2O^{-1}\partial\mathbf{w}\big)\nm\\
  &=& \lambda D\partial\mathbf{w}   + \lambda\big( OK\mathrm{U}\mathrm{V}K(\mathbf{U}-\mathrm{1})+OK\mathrm{U}K\mathrm{U}O^{-1}\partial\mathbf{w}-D\partial\mathbf{w}\big).\label{n4}
 \eer

 Then as  $|x|$ is sufficiently large we   get
 \ber
&&\Delta |(\partial\mathbf{w})|^2\nm\\
&&\ge2(\partial\mathbf{w})^\tau\Delta (\partial\mathbf{w})\nm\\
&&= \lambda (\partial\mathbf{w})^\tau D\partial\mathbf{w}+ \lambda(\partial\mathbf{w})^\tau\big( OK\mathrm{U}\mathrm{V}K(\mathbf{U}-\mathrm{1})+OK\mathrm{U}K\mathrm{U}O^{-1}\partial\mathbf{w}-D\partial\mathbf{w}\big).\label{a43}
 \eer
Hence, similar to  \eqref{a40},  for any sufficiently small $\vep\in(0, 1)$, there exists an $R_\vep>R$ such that
\ber
   \Delta |\partial\mathbf{w}|^2\ge 2\lambda\sigma_0^2\left(1-\frac\vep2\right)|\partial\mathbf{w}|^2, \quad |x|>R_\vep. \label{a44}
 \eer
Since we have shown that $|\partial\mathbf{w}|^2\to 0$ as $|x|\to \infty$,  from \eqref{a44} we get that there exists a positive constant $C(\vep)$ such that
 \ber
   |\partial\mathbf{w}|^2\le C(\vep)\re^{-\sigma_0\sqrt{2\lambda}(1-\vep)|x|}, \quad |x|>R_\vep,\label{a45}
  \eer
which gives \eqref{d2}.

Using the decay estimates   one can now calculate the magnetic fluxes.
Indeed,  from  the equations \eqref{F1}--\eqref{F2} one has
\ber
F_{12}^0 &=& -\frac{1}{\sqrt{2N}} \Delta\left([N-1]\ln|\phi|^2 +|\ln\phi_N|^2\right),\\
F_{12}^{N^2-1} &=& -\sqrt{\frac{N-1}{2N}} \Delta\left(\ln|\phi|^2 - \ln\phi_N|^2\right),
 \eer
so that
  \ber
   {\cal F}^{U(1)}  &=&-\frac{1}{\sqrt{2N}} \itr\Delta \left([N-1]\ln|\phi|^2+|\ln\phi_N|^2\right)\ud x,\label{d5}\\
     {\cal F}^{SU(N)}  &=&-\sqrt{\frac{N-1}{2N}} \itr\Delta \left(\ln|\phi|^2-|\ln\phi_N|^2\right)\ud x, \label{d5a}
  \eer
and direct integration leads to   \eqref{d6}--\eqref{d6a} which show that the vortex magnetic fluxes are completely determined by the zeros of the Higgs scalar.

 Then we complete the  proof of Theorem \ref{th2} and Theorem \ref{th3} for the planar case.

 \section{Doubly periodic case}
 \setcounter{equation}{0}
  \setcounter{lemma}{0}

In this section we establish the existence of doubly periodic solutions for \eqref{a1}--\eqref{a2}. We will use a constrained minimization approach, developed in \cite{CaffarelliYang} and later refined by
\cite{taran96,nota}, to establish the existence of the first solution.  The key step is to find out an inequality type of constraints.     We show that when the coupling parameter $\lambda$ is sufficiently large the variational problem with such constraints admits an  interior critical point, which is also a critical point of the original variational problem. To get the existence of the second solution, we use the mountain pass theorem.

Consider the equations \eqref{a1}--\eqref{a2} over a doubly periodic domain $\Omega$.  We first derive  a priori estimate for the solutions to \eqref{a1}--\eqref{a2}.
 \begin{lemma}\label{lem1}
  Every solution $(u_1, u_2)$ to \eqref{a1}--\eqref{a2} satisfies
  \[u_1<0, \quad u_2<0    \quad \text{in} \quad \Omega.\]
 \end{lemma}

{\bf Proof.} \quad  To prove this lemma, it is convenient to rewrite \eqref{a1}--\eqref{a2} equivalently as
 \ber
  \Delta u_1&=&\lambda\left\{\frac{\kappa-1}{N^2}\big([N-1+\kappa]\re^{u_1}+[N-1][\kappa-1]\re^{u_2}\big)\big(\re^{u_1}-\re^{u_2}\big)\right.\nm\\
  &&\left.+\frac1N\big([N-1+\kappa]\re^{u_1}[\re^{u_1}-1]+[\kappa-1]\re^{u_2}[1-\re^{u_2}]\big)\right\} +4\pi\sum\limits_{s=1}^{n_1}\delta_{p_{1s}}, \label{3a}\\
  \Delta u_2&=&\lambda\left\{\frac{(N-1)(\kappa-1)}{N^2}\big([\kappa-1]\re^{u_1}+[1+[N-1]\kappa]\re^{u_2}\big)\big(\re^{u_2}-\re^{u_1}\big)\right.\nm\\
  &&\left.+\frac1N\big([1+[N-1]\kappa]\re^{u_2}[\re^{u_2}-1]+[N-1][\kappa-1]\re^{u_1}[1-\re^{u_1}]\big)\right\} +4\pi\sum\limits_{s=1}^{n_2}\delta_{p_{2s}}.\label{3b}
 \eer
 It is easy to see that $u_i$  may achieve its maximum value at some point  $ \tilde{x}_i\in\Omega\setminus\{p_{i1}, \dots, p_{in_i}\},\, i=1,2$.
 Denote $\tilde{u}_i\equiv \max\limits_{x\in \Omega}u_i=u_i(\tilde{x}_i), \, i=1,2.$

  We first show that $\tilde{u}_i\le 0, \, i=1, 2.$  For the  case $\tilde{u}_1\ge \tilde{u}_2$,  using \eqref{3a}, we have
   \berr
    0\ge \Delta u_1(\tilde{x}_1)&=&\lambda\left\{\frac{\kappa-1}{N^2}\left([N-1+\kappa]\re^{\tilde{u}_1}+[N-1][\kappa-1]\re^{u_2(\tilde{x}_1)}\right)\left(\re^{\tilde{u}_1}-\re^{u_2(\tilde{x}_1)}\right)\right.\nm\\
  &&\left.+\frac1N\left([N-1+\kappa]\re^{\tilde{u}_1}\left[\re^{\tilde{u}_1}-1\right]+[\kappa-1]\re^{u_2(\tilde{x}_1)}\left[1-\re^{u_2(\tilde{x}_1)}\right]\right)\right\} \nm\\
   &\ge&  \frac{\lambda}{N}\left( [N-1+\kappa]\re^{\tilde{u}_1}- [\kappa-1]\re^{u_2(\tilde{x}_1)}\right)\left(\re^{\tilde{u}_1}-1\right).
   \eerr
Then from maximum principle we see that $\tilde{u}_1\le0$.   Then in this case both  $\tilde{u}_1$ and  $\tilde{u}_2$ are nonpositive.

If   $\tilde{u}_2\ge \tilde{u}_1$,  using \eqref{3b}, we obtain
    \berr
    0\ge \Delta u_2(\tilde{x}_2)&=&\lambda\left\{\frac{(N-1)(\kappa-1)}{N^2}\left([\kappa-1]\re^{u_1(\tilde{x}_2)}+[1+[N-1]\kappa]\re^{\tilde{u}_2}\right)\left(\re^{\tilde{u}_2}-\re^{u_1(\tilde{x}_2)}\right)\right.\nm\\
   &&\left.+\frac1N\left([1+[N-1]\kappa]\re^{\tilde{u}_2}\left[\re^{\tilde{u}_2}-1\right]+[N-1][\kappa-1]\re^{u_1(\tilde{x}_2)}\left[1-\re^{u_1(\tilde{x}_2)}\right]\right)\right\}\\
   &\ge&  \frac{\lambda}{N}\left( [1+[N-1]\kappa]\re^{\tilde{u}_2}-[N-1][\kappa-1]\re^{u_1(\tilde{x}_2)}\right)\left(\re^{\tilde{u}_2}-1\right).
   \eerr
Using maximum principle again, we get $\tilde{u}_2\le 0$. Hence  the  conclusion follows for this case. Therefore we have $u_i\le 0, \, i=1,2.$

To prove the strict inequality, we can  apply the strong maximum principle. In fact,  it is sufficient to note that
 \berr
  \Delta u_1+a_1(x)u_1&=&\frac{\lambda}{N^2}(\kappa-1)\left( [N-1+\kappa]\re^{u_1}+(1+[N-1]\kappa)\re^{u_2} \right)(1-\re^{u_2})\ge0,\\
   \Delta u_2+a_2(x)u_2&=&\frac{\lambda}{N^2}(N-1)(\kappa-1)\left( [N-1+\kappa]\re^{u_1}+(1+[N-1]\kappa)\re^{u_2} \right)(1-\re^{u_1})\ge0,
 \eerr
 where
 \berr
   a_1(x)&=&\frac{\lambda}{N^2}\left([N-1+\kappa]^2\re^{u_1}+[N-1][\kappa-1]^2\re^{u_2}\right)\frac{1-\re^{u_1}}{u_1},\\
    a_2(x)&=&\frac{\lambda}{N^2}\left([N-1][\kappa-1]^2\re^{u_1}+[1+[N-1][\kappa-1]]^2\re^{u_2}\right)\frac{1-\re^{u_2}}{u_2}.
 \eerr
 Then  Lemma \ref{lem1} follows from the strong maximum principle.

From Lemma \ref{lem1}, we get the first part of Theorem \ref{th2}.

Let $u_i^0$ be the solution of the following problem (see \cite{aubi})
\berr
 &&\Delta u_i^0=4\pi\sum\limits_{s=1}^{n_i}\delta_{p_{is}}-\frac{4\pi n_i}{|\Omega|},\\
  &&\ito u_i^0\ud x=0, \quad i=1, 2,
\eerr
and $u_i=u_i^0+v_i, \, i=1, 2$. Then  over $\Omega$ the equations \eqref{a1}--\eqref{a2} can be reduced as
  \ber \Delta v_1&=&\lambda\left\{\frac{1}{N^2}\Big([N-1+\kappa]^2\re^{2u_1^0+2v_1}-[\kappa-1]\big(N-[N-2][\kappa-1]\big)\re^{u_1^0+u_2^0+v_1+v_2}\right.\nm\\
 &&\left.+[1-\kappa]\big(1+[N-1]\kappa\big)\re^{2u_2^0+2v_2}\Big)-\frac1N\Big([N-1+\kappa]\re^{u_1^0+v_1}+[1-\kappa]\re^{u_2^0+v_2}\Big)\right\}\nm\\
 &&+\frac{4\pi n_1}{|\Omega|},\label{c1}\\
 \Delta v_2&=&\lambda\left\{\frac{1}{N^2}\Big([N-1][1-\kappa][N-1+\kappa]\re^{2u_1^0+2v_1}-[N-1][\kappa-1]\big(2+[N-2]\kappa\big)\re^{u_1^0+u_2^0+v_1+v_2}\right.\nm\\
&&\left.
 +\big(1+[N-1]\kappa\big)^2\re^{2u_2^0+2v_2}\Big)-\frac1N\Big([N-1][1-\kappa]\re^{u_1^0+v_1}+\big(1+[N-1]\kappa\big)\re^{u_2^0+v_2}\Big)\right\}\nm\\
   &&+\frac{4\pi n_2}{|\Omega|}.\label{c2}
\eer

To find a variational principle for the problem \eqref{c1}--\eqref{c2}, as in the full plane case we rewrite the equations \eqref{c1}--\eqref{c2} equivalently as
\ber
 \left(N-1+\frac1\kappa\right)\Delta v_1+\left(1-\frac1\kappa\right)\Delta v_2&=&\lambda\left([N-1+\kappa]\re^{2u_1^0+2v_1}-N\re^{u_1^0+v_1}\right.\nm\\
 &&\left.-[\kappa-1]\re^{u_1^0+u_2^0+v_1+v_2}\right)+ \frac{b_1}{|\Omega|},\label{c3}\\
\left(1-\frac1\kappa\right)\Delta v_1+\left(\frac{1}{N-1}+\frac1\kappa\right)\Delta v_2&=&\lambda\left(\left[\frac{1}{N-1}+\kappa\right]\re^{2u_2^0+2v_2}-\frac{N}{N-1}\re^{u_2^0+v_2}
\right.\nm\\
 &&\left.-[\kappa-1]\re^{u_1^0+u_2^0+v_1+v_2}\right)+ \frac{b_2}{|\Omega|},\label{c4}
\eer
where the notation
\ber
 b_1\equiv\frac{4\pi\big([1+[N-1]\kappa]n_1+[\kappa-1]n_2\big)}{\kappa}, \quad b_2\equiv\frac{4\pi\big([N-1][\kappa-1]n_1+[N-1+\kappa]n_2\big)}{(N-1)\kappa} \label{c4'}
\eer
will be used throughout this paper.

We will work on the space $W^{1, 2}(\Omega)\times W^{1, 2}(\Omega)$, where $W^{1, 2}(\Omega)$ is  the set of $\Omega$-periodic $L^2$-
functions whose  derivatives  also belong  $L^2(\Omega)$. We denote the usual  norm on $W^{1, 2}(\Omega)$  by $\|\cdot\|$ as given by
 $\|w\|^2=\|w\|_2^2+\|\nabla w\|_2^2=\ito w^2\ud x+\ito |\nabla w|^2\ud x$.

Then we  easily see  that the equations \eqref{c3}--\eqref{c4} are the Euler--Lagrange equations of the functional
  \ber
   I(v_1, v_2)&=&\frac12\left(N-1+\frac1\kappa\right)\|\nabla v_1\|_2^2+\frac12\left(\frac{1}{N-1}+\frac1\kappa\right)\|\nabla v_2\|_2^2+\left(1-\frac1\kappa\right)\ito\nabla v_1\cdot\nabla v_2\ud x\nm\\
   &&+\frac\lambda2\left([N-1+\kappa]\ito[\re^{u_1^0+v_1}-1]^2\ud x+\left[\frac{1}{N-1}+\kappa\right]\ito[\re^{u_2^0+v_2}-1]^2\ud x\right.\nm\\
   &&\left.+2[1-\kappa]\ito[\re^{u_1^0+v_1}-1][\re^{u_2^0+v_2}-1]\ud x\right)+\frac{b_1}{|\Omega|}\ito v_1\ud x+\frac{b_2}{|\Omega|}\ito v_2\ud x. \label{c5}
  \eer
Hence in the following subsections we concentrate on finding the critical points of the functional $I$.

\subsection{Constrained minimization procedure}
To find a first critical point  of the functional $I$, we carry out a constrained minimization procedure.

For any solution $(v_1, v_2)$ of \eqref{c3}--\eqref{c4}, integrating over $\Omega$  gives the following constraints
\ber
(N-1+\kappa)\ito\re^{2u_1^0+2v_1}\ud x-N\ito\re^{u_1^0+v_1}\ud x-(\kappa-1)\ito\re^{u_1^0+u_2^0+v_1+v_2}\ud x+\frac{b_1}{\lambda}=0,\label{c6}\\
 \left(\frac{1}{N-1}+\kappa\right)\ito\re^{2u_2^0+2v_2}-\frac{N}{N-1}\ito\re^{u_2^0+v_2}-(\kappa-1)\ito\re^{u_1^0+u_2^0+v_1+v_2}\ud x+\frac{b_2}{\lambda}=0.\label{c7}
\eer

We first establish the  necessary condition   stated in Theorem \ref{th2} for the existence of  solutions to \eqref{a1}--\eqref{a2}.  To this end,  for any solution $(v_1, v_2)$ of \eqref{c3}--\eqref{c4},   from  \eqref{c6}--\eqref{c7}  we observe that
\ber
 &&\ito \left(\left[N-1+\kappa\right]\re^{2u_1^0+2v_1}+\left[\frac{1}{N-1}+\kappa\right]\re^{2u_2^0+2v_2}-2[\kappa-1]\re^{u_1^0+u_2^0+v_1+v_2}\right.\nm\\
 &&\left.\quad-N\re^{u_1^0+v_1}-\frac{N}{N-1}\re^{u_2^0+v_2}\right)\ud x =-\frac{4\pi N([N-1]n_1+n_2)}{(N-1)\lambda},\label{c7f}
\eer
where we used the notation \eqref{c4'}.

 Consider a function
 \be
 g(t_1, t_2)=(N-1+\kappa)t_1^2+\left(\frac{1}{N-1}+\kappa\right)t_2^2-2(\kappa-1)t_1t_2-Nt_1-\frac{N}{N-1}t_2, (t_1, t_2)\in \mathbb{R}^2.\label{c7h}
 \ee
 Notice that the Hessian of $g$ is $\frac{2}{N} A(N,\frac{1}\kappa)$ where the matrix $A (N,\kappa)$ was defined in  \ref{11}. We may check that the function  $g$ reaches its unique  global minimum
 \be g_{\min}=-\frac{N^2}{4(N-1)}  \quad   \text{at} \quad  \left(\frac12,\, \frac12\right). \label{c7g}
 \ee
Then  we conclude from  \eqref{c7f}--\eqref{c7g} that
 \[-\frac{N^2|\Omega|}{4(N-1)}=\ito g_{\min}\ud x\le -\frac{4\pi N([N-1]n_1+n_2)}{(N-1)\lambda}, \]
 which implies  the Bradlow's bound \cite{masu}
 \[\lambda\ge \frac{16\pi([N-1]n_1+n_2)}{N|\Omega|}. \]
Then we get the necessary condition  \eqref{a9'}  stated  in Theorem \ref{th2}.

Using \eqref{c6}--\eqref{c7}, we  also find
 \ber
  && \frac12\left([N-1+\kappa]\ito\left[\re^{u_1^0+v_1}-1\right]^2\ud x+\left[\frac{1}{N-1}+\kappa\right]\ito\left[\re^{u_2^0+v_2}-1\right]^2\ud x\right.\nm\\
   &&\left.+2[1-\kappa]\ito\left[\re^{u_1^0+v_1}-1\right]\left[\re^{u_2^0+v_2}-1\right]\ud x\right)\nm\\
   &&=\frac12\left(N\ito\left[1-\re^{u_1^0+v_1}\right]\ud x+\frac{N}{N-1}\ito\left[1-\re^{u_2^0+v_2}\right]\ud x\right)-\frac{2\pi N}{\lambda}\left(n_1+\frac{n_2}{N-1}\right). \label{c8}
   \eer

It is well-known that the space  $W^{1,2}(\Omega)$ can be decomposed  as
  \[
  W^{1,2}(\Omega)=\mathbb{R}\oplus  \dot{W}^{1,2}(\Omega).
  \]
  where
\[\dot{W}^{1,2}(\Omega)=\left\{w\in W^{1,2}(\Omega)\Bigg| \ito w\ud x=0\right\}\]
  is a closed subspace of $W^{1,2}(\Omega)$.

Then, for $v_i\in W^{1,2}(\Omega)$ we  have the decomposition
 \[v_i=c_i+w_i, \]
where
  \[\ito w_i\ud x=0,\quad c_i=\frac{1}{|\Omega|}\ito v_i \ud x,\quad i=1,2. \]

If $(v_1, v_2)\in W^{1,2}(\Omega)\times W^{1,2}(\Omega)$ satisfies  \eqref{c6}--\eqref{c7}, we obtain
 \ber
  (N-1+\kappa)\re^{2c_1}\ito\re^{2u_1^0+2w_1}\ud x-\re^{c_1}P_1(w_1, w_2, \re^{c_2})+\frac{b_1}{\lambda}=0,\label{c9}\\
  \left(\frac{1}{N-1}+\kappa\right)\re^{2c_2}\ito\re^{2u_2^0+2w_2}\ud x-\re^{c_2}P_2(w_1, w_2, \re^{c_1})+\frac{b_2}{\lambda}=0,\label{c10}
 \eer
where
\ber
P_1(w_1, w_2, \re^{c_2})&\equiv& N\ito\re^{u_1^0+w_1}\ud x+(\kappa-1)\re^{c_2}\ito\re^{u_1^0+u_2^0+w_1+w_2}\ud x,\label{c11}\\
P_2(w_1, w_2, \re^{c_1})&\equiv& \frac{N}{N-1}\ito\re^{u_2^0+w_2}+(\kappa-1)\re^{c_1}\ito\re^{u_1^0+u_2^0+w_1+w_2}\ud x.\label{c12}
\eer

Then \eqref{c9}--\eqref{c10} are solvable with respect to $c_1, c_2$ if and only if
 \ber
   P_1^2(w_1, w_2, \re^{c_2})&\ge& \frac{4b_1}{\lambda}(N-1+\kappa)\ito\re^{2u_1^0+2w_1}\ud x, \label{c13}\\
   P_2^2(w_1, w_2, \re^{c_1})&\ge&\frac{4b_2}{\lambda}\left(\frac{1}{N-1}+\kappa\right)\ito\re^{2u_2^0+2w_2}\ud x. \label{c14}
 \eer

In view of \eqref{c13}--\eqref{c14}, we choose the following inequality type constraints
 \ber
   \left(\ito\re^{u_1^0+w_1}\ud x\right)^2&\ge& \frac{4(N-1+\kappa)b_1}{N^2\lambda}\ito\re^{2u_1^0+2w_1}\ud x,\label{c16}\\
   \left(\ito\re^{u_2^0+w_2}\ud x\right)^2&\ge&  \frac{4(N-1)(1+[N-1]\kappa)b_2}{N^2\lambda}\ito\re^{2u_2^0+2w_2}\ud x.\label{c17}
 \eer

We introduce the admissible set as
 \be
  \mathcal{A}=\Big\{(w_1, w_2)\in \dot{W}^{1, 2}(\Omega)\times \dot{W}^{1, 2}(\Omega)\, \,\text{satisfies}\,\, \eqref{c16},\, \eqref{c17}\Big\}. \label{c18}
 \ee

Then,  for any  $(w_1, w_2)\in\mathcal{A}$, we can obtain a solution of  \eqref{c9}--\eqref{c10} with respect to $c_1$ and $c_2$ by solving    the following equations
 \ber
 \re^{c_1}&=&\frac{ P_1(w_1, w_2,  \re^{c_2})+\sqrt{P_1^2(w_1,  w_2,  \re^{c_2})-\frac{4(N-1+\kappa)b_1}{\lambda}\ito\re^{2u_1^0+2w_1}\ud x}}{2(N-1+\kappa)\ito\re^{2u_1^0+2w_1}\ud x}
 \nm\\&\equiv& f_1(\re^{c_2}), \label{c19}\\
  \re^{c_2}&=&\frac{ P_2(w_1, w_2,  \re^{c_1})+\sqrt{P_2^2(w_1, w_2,  \re^{c_1})-\frac{4(1+[N-1]\kappa)b_2}{(N-1)\lambda}\ito\re^{2u_2^0+2w_2}\ud x}}{2\left(\frac{1}{N-1}+\kappa\right)\ito\re^{2u_2^0+2w_2}\ud x}
 \nm\\&\equiv& f_2(\re^{c_1}).\label{c20}
 \eer

 To  solve \eqref{c19}--\eqref{c20},  we just need to   find  the zeros of the function
   \[ f(X)\equiv X-f_1(f_2(X)), \quad X\ge0.  \]
 In order to do that,  it sufficient to prove the following proposition.
 \begin{proposition}\label{prop1}
 For any $(w_1, w_2)\in   \mathcal{A}$, the equation
  \[f(X)=X-f_1(f_2(X))=0\]
   admits a unique positive solution $X_0$.
 \end{proposition}

For any  $(w_1, w_2)\in\mathcal{A}$,  using  this proposition,    we can get a solution of  \eqref{c9}--\eqref{c10} with respect to $ c_1, c_2$.

 {\bf Proof of Proposition \ref{prop1}. }   \quad  We see from \eqref{c19} and \eqref{c20} that
 \be
 f_i(X)>0, \quad \forall\, X\ge0, \, i=1,2.\label{c21}
 \ee
 Hence  $f(0)=-f_1(f_2(0))<0$.  It is easy to check  that
 \ber
  \frac{\ud f_1(X)}{\ud X}&=& \frac{(\kappa-1) f_1(X)\ito\re^{u_1^0+u_2^0+w_1+w_2}\ud x}{\sqrt{P_1^2(w_1, w_2,  X)-\frac{4(N-1+\kappa)b_1}{\lambda}\ito\re^{2u_1^0+2w_1}\ud x}}, \label{c22}\\
  \frac{\ud f_2(X)}{\ud X}&=& \frac{(\kappa-1) f_2(X)\ito\re^{u_1^0+u_2^0+w_1+w_2}\ud x}{\sqrt{P_2^2(w_1, w_2,  X)-\frac{4(1+[N-1]\kappa)b_2}{(N-1)\lambda}\ito\re^{2u_2^0+2w_2}\ud x}},\label{c23}
 \eer
  which are all positive since $\kappa>1$.  Thus,  the functions $f_i(X)$, $(i=1, 2)$ are  strictly increasing for  all  $X>0$.

  A direct computation gives
    \berr
    \lim\limits_{X\to+\infty}\frac{f_1(X)}{X}&=&\frac{(\kappa-1)\ito\re^{u_1^0+u_2^0+w_1+w_2}\ud x}{(N-1+\kappa)\ito\re^{2u_1^0+2w_1}\ud x}, \\
    \lim\limits_{X\to+\infty}\frac{f_2(X)}{X}&=&\frac{(\kappa-1)\ito\re^{u_1^0+u_2^0+w_1+w_2}\ud x}{\left(\frac{1}{N-1}+\kappa\right)\ito\re^{2u_2^0+2w_2}\ud x},
    \eerr
which implies
   \berr
  \lim\limits_{X\to+\infty}\frac{f(X)}{X}&=&1-\frac{(\kappa-1)^2\left(\ito\re^{u_1^0+u_2^0+w_1+w_2}\ud x\right)^2}{(N-1+\kappa)\left(\frac{1}{N-1}+\kappa\right)\ito\re^{2u_1^0+2w_1}\ud x\ito\re^{2u_2^0+2w_2}\ud x}\\
   &\ge& 1-\frac{(\kappa-1)^2}{(N-1+\kappa)\left(\frac{1}{N-1}+\kappa\right)}\\
   &=&\frac{N^2\kappa}{(N-1)(\kappa-1)^2+N^2\kappa}>0
   \eerr
 Hence we conclude that  \[\lim\limits_{X\to+\infty}f(X)=+\infty.\]

In view of the fact  $f(0)<0$,  then  we  infer that the function  $f(\cdot)$ has at least one zero point $X_0>0$.

 Now we  prove that the  zero point of $f(\cdot)$ is also unique. From \eqref{c22}--\eqref{c23} and \eqref{c19}--\eqref{c20} we obtain
 \berr
  \frac{\ud f(X)}{\ud X}&=&1-\frac{(\kappa-1)^2f_1(f_2(X))f_2(X)\left(\ito\re^{u_1^0+u_2^0+w_1+w_2}\ud x\right)^2}
  {\sqrt{P_1^2(w_1, w_2,  X)-\frac{4(N-1+\kappa)b_1}{\lambda}\ito\re^{2u_1^0+2w_1}\ud x}}\times
\\&&\times\frac{1}{\sqrt{P_2^2(w_1, w_2,  X)-\frac{4(1+[N-1]\kappa)b_2}{(N-1)\lambda}\ito\re^{2u_2^0+2w_2}\ud x}}\\
  &>&1-\frac{f_1(f_2(X))}{X}=\frac{f(X)}{X},
 \eerr
which  gives    \[  \left(\frac{f(X)}{X}\right)'>0. \]
  Thus we  see that  $\frac{f(X)}{X}$ is strictly increasing  for $X>0$.  Consequently, $f(X)$ in strictly increasing for $X>0$, which implies  $f(X)$ has a unique zero point.
 Then the proof of Proposition \ref{prop1} is complete.

 From the above discussion we conclude  that, for any $(w_1, w_2)\in \mathcal{A}$,   there is  a  pair $(c_1(w_1, w_2), c_2(w_1, w_2))$ given by
 \eqref{c19}--\eqref{c20},  which solves \eqref{c9}--\eqref{c10},  such that  $(v_1, v_2)$ defined by
 \[ v_i=w_i+c_i(w_1, w_2),  \quad i=1,2\]
 satisfies  \eqref{c6}--\eqref{c7}.

Thus, to seek the critical points of $I$, we may consider the functional
  \be
   J(w_1, w_2)=I(w_1+c_1(w_1, w_2), w_2+c_2(w_1, w_2)), \quad  (w_1, w_2)\in \mathcal{A}. \label{c24}
  \ee

In view of \eqref{c8}, we may write the functional $J$ as
 \ber
  J(w_1, w_2)&=&\frac12\left(N-1+\frac1\kappa\right)\|\nabla w_1\|_2^2+\frac12\left(\frac{1}{N-1}+\frac1\kappa\right)\|\nabla w_2\|_2^2+\left(1-\frac1\kappa\right)\ito\nabla w_1\cdot\nabla w_2\ud x\nm\\
   &&+\frac\lambda2\left(N\ito\left[1-\re^{c_1}\re^{u_1^0+w_1}\right]\ud x+\frac{N}{N-1}\ito\left[1-\re^{c_2}\re^{u_2^0+w_2}\right]\ud x\right)\nm\\
   &&-2\pi N\left(n_1+\frac{n_2}{N-1}\right)+b_1c_1+b_2c_2, \label{c24'}
 \eer
where $b_1, b_2$ are defined by \eqref{c4'}.

We easily see that the functional $J$ is Frech\'{e}t differentiable in the interior of $\mathcal{A}$.  If we find a critical point $(w_1, w_2)$ of
$J$, which lies in the interior of $\mathcal{A}$, then $(w_1+c_1(w_1, w_2), w_2+c_2(w_1, w_2))$ is  a critical point of $I$.  Therefore in what follows  we
 just need to find  the critical points for the functional  $J$.

 We  first establish  the following lemma.
  \begin{lemma}\label{lem2}
   For any $(w_1, w_2)\in \mathcal{A}$,  there holds
    \ber
     \re^{c_i}\ito\re^{u_i^0+w_i}\le |\Omega|, \quad i=1, 2.\label{c25}\\
      \re^{c_i}\le 1, \quad i=1, 2. \label{c25'}
    \eer
  \end{lemma}

 {\bf Proof.} \quad For any  $(w_1, w_2)\in \mathcal{A}$,  from \eqref{c19}--\eqref{c20} we  obtain
  \ber
   \re^{c_1}&\le& \frac{N\ito\re^{u_1^0+w_1}\ud x+(\kappa-1)\re^{c_2}\ito\re^{u_1^0+u_2^0+w_1+w_2}\ud x}{(N-1+\kappa)\ito\re^{2u_1^0+2w_1}\ud x},\label{c26}\\
   \re^{c_2}&\le&\frac{N\ito\re^{u_2^0+w_2}\ud x+(N-1)(\kappa-1)\re^{c_1}\ito\re^{u_1^0+u_2^0+w_1+w_2}\ud x}{\left(1+[N-1]\kappa\right)\ito\re^{2u_2^0+2w_2}\ud x}.\label{c27}
  \eer
  Then, using \eqref{c26}--\eqref{c27} and H\"{o}lder inequality we have
   \berr
     \re^{c_1}&\le& \frac{N\ito\re^{u_1^0+w_1}\ud x}{{(N-1+\kappa)\ito\re^{2u_1^0+2w_1}\ud x}}\nm\\
     &&+\frac{N(\kappa-1) \ito\re^{u_2^0+w_2}\ud x\ito\re^{u_1^0+u_2^0+w_1+w_2}\ud x}{(N-1+\kappa)\left(1+[N-1]\kappa\right)\ito\re^{2u_1^0+2w_1}\ud x\ito\re^{2u_2^0+2w_2}\ud x}\nm\\
     &&+\frac{(N-1)(\kappa-1)^2\re^{c_1}\left(\ito\re^{u_1^0+u_2^0+w_1+w_2}\ud x\right)^2}{(N-1+\kappa)\left(1+[N-1]\kappa\right)\ito\re^{2u_1^0+2w_1}\ud x\ito\re^{2u_2^0+2w_2}\ud x}\nm\\
     &\le& \frac{N\ito\re^{u_1^0+w_1}\ud x}{{(N-1+\kappa)\ito\re^{2u_1^0+2w_1}\ud x}}+\frac{(N-1)(\kappa-1)^2\re^{c_1}}{(N-1+\kappa)\left(1+[N-1]\kappa\right)}\nm\\
     &&+\frac{N(\kappa-1) \ito\re^{u_2^0+w_2}\ud x\ito\re^{u_1^0+u_2^0+w_1+w_2}\ud x}{(N-1+\kappa)\left(1+[N-1]\kappa\right)\ito\re^{2u_1^0+2w_1}\ud x\ito\re^{2u_2^0+2w_2}\ud x},
   \eerr
which  implies
 \ber
  \re^{c_1}&\le&\frac{\left(1+[N-1]\kappa\right)\ito\re^{u_1^0+w_1}\ud x}{N\kappa\ito\re^{2u_1^0+2w_1}\ud x}+\frac{(\kappa-1)\ito\re^{u_2^0+w_2}\ud x\ito\re^{u_1^0+u_2^0+w_1+w_2}\ud x}{N\kappa\ito\re^{2u_1^0+2w_1}\ud x\ito\re^{2u_2^0+2w_2}\ud x}.\label{c28}
 \eer
Analogously,  we
 \ber
  \re^{c_2}&\le&\frac{\left(N-1+\kappa\right)\ito\re^{u_2^0+w_2}\ud x}{N\kappa\ito\re^{2u_2^0+2w_2}\ud x}+\frac{(N-1)(\kappa-1)\ito\re^{u_1^0+w_1}\ud x\ito\re^{u_1^0+u_2^0+w_1+w_2}\ud x}{N\kappa\ito\re^{2u_1^0+2w_1}\ud x\ito\re^{2u_2^0+2w_2}\ud x}.\label{c29}
 \eer
   Then   from \eqref{c28} and the H\"{o}lder inequality we see that
    \berr
       \re^{c_1}\ito\re^{u_1^0+w_1}\ud x&\le&  \frac{\left(1+[N-1]\kappa\right)\left(\ito\re^{u_1^0+w_1}\ud x\right)^2}{N\kappa\ito\re^{2u_1^0+2w_1}\ud x}
      \nm\\
      &&+\frac{(\kappa-1)\ito\re^{u_1^0+w_1}\ud x\ito\re^{u_2^0+w_2}\ud x\ito\re^{u_1^0+u_2^0+w_1+w_2}\ud x}{N\kappa\ito\re^{2u_1^0+2w_1}\ud x\ito\re^{2u_2^0+2w_2}\ud x}\nm\\
      &\le&\frac{\left(1+[N-1]\kappa+\kappa-1\right)|\Omega|}{N\kappa}=|\Omega|.
    \eerr

Analogously,  we  get
  \berr
       \re^{c_2}\ito\re^{u_2^0+w_2}\ud x&\le& |\Omega|.
  \eerr
Then   we obtain  \eqref{c25}, which implies \eqref{c25'}  by using  Jensen's  inequality.  The proof of Lemma \ref{lem2} is complete.

\begin{lemma}\label{lem3}
 For any $(w_1,w_2)\in\mathcal{A}$  and  $s\in (0, 1)$, there holds
  \ber
   \ito\re^{u_1^0+w_1}\ud x&\le& \left(\frac{N^2\lambda}{4[N-1+\kappa]b_1}\right)^{\frac{1-s}{s}}\left(\ito\re^{su_1^0+sw_1}\ud x\right)^{\frac1s},\label{c30}\\
    \ito\re^{u_2^0+w_2}\ud x&\le& \left(\frac{N^2\lambda}{4[N-1](1+[N-1]\kappa)b_2}\right)^{\frac{1-s}{s}}\left(\ito\re^{su_2^0+sw_2}\ud x\right)^{\frac1s}.\label{c31}
  \eer
\end{lemma}

{\bf Proof.}   To prove this lemma,  we use the approach developed in \cite{nota, nota2}.

For $s\in (0, 1)$, let $\gamma=\frac{1}{2-s}$ such that $   s\gamma+2(1-\gamma)=1$. In view of the H\"{o}lder inequality, we have
 \berr
  \ito\re^{u_1^0+w_1}\ud x&\le&\left(\ito\re^{su_1^0+sw_1}\ud x\right)^\gamma\left(\ito\re^{2u_1^0+2w_1}\ud x\right)^{1-\gamma} \\
  &\le& \left(\frac{N^2\lambda}{4[N-1+\kappa]b_1}\right)^{1-\gamma}\left(\ito\re^{u_1^0+w_1} \ud x\right)^{2(1-\gamma)}\left(\ito\re^{su_1^0+sw_1} \ud x\right)^\gamma,
 \eerr
which gives
 \berr
   \ito\re^{u_1^0+w_1}\ud x&\le&  \left(\frac{N^2\lambda}{4[N-1+\kappa]b_1}\right)^{\frac{1-\gamma}{2\gamma-1}}\left(\ito\re^{su_1^0+sw_1} \ud x\right)^{\frac{\gamma}{2\gamma-1}}\\
   &=& \left(\frac{N^2\lambda}{4[N-1+\kappa]b_1}\right)^{\frac{1-s}{s}}\left(\ito\re^{su_1^0+sw_1} \ud x\right)^{\frac{1}{s}}.
 \eerr
Then \eqref{c30} is established. Similarly, we can prove \eqref{c31}. The proof of Lemma \ref{lem3} is complete.

 Next, we use Lemma \ref{lem3} to   show that the functional $J$ is coercive on $ \mathcal{A}$. To this end, we need the  Moser--Trudinger inequality
(see\cite{aubi,font})
 \be
 \ito \re^{w}\ud x \le C_1\exp\left(\frac{1}{16\pi}\|\nabla w\|_2^2\right), \quad \forall\, w\in \dot{W}^{1, 2}(\Omega),\label{c32}
 \ee
 where $C_1$ is a positive constant depending on $\Omega$ only.

\begin{lemma}\label{lem4}
  For any $(w_1, w_2)\in \mathcal{A}$, there exist suitable positive constants $C_2$ and $C_3$, independent of $\lambda$,  such that
  \ber
   J(w_1, w_2)\ge C_2(\|\nabla w_1\|_2^2+\|\nabla w_2\|_2^2)-C_3(\ln\lambda+1). \label{c32'}
  \eer
\end{lemma}

{\bf Proof.}\quad
 Noting that the matrices $A(N, \kappa)$ and  $A(N, \kappa^{-1})$  defined by \eqref{11} are both  positive definite,
then  using  \eqref{c8}, \eqref{c24'}, we have
 \ber
  J(w_1, w_2)&\ge&\frac{\alpha_0(\kappa)}{2}\left(\|\nabla w_1\|_2^2+\|\nabla w_2\|_2^2\right)+b_1c_1+b_2c_2, \label{c33}
 \eer
  where $\alpha_0(\kappa)$ is a positive constant defined by \eqref{a140}.

Next we estimate $c_1, c_2$ in \eqref{c33}.  We see from \eqref{c19}--\eqref{c20} that
\berr
 \re^{c_1}&\ge& \frac{N\ito\re^{u_1^0+w_1}\ud x}{2(N-1+\kappa)\ito\re^{2u_1^0+2w_1}\ud x},\\
 \re^{c_2}&\ge& \frac{N\ito\re^{u_2^0+w_2}\ud x}{2(1+[N-1]\kappa)\ito\re^{2u_2^0+2w_2}\ud x}.
\eerr
Then,  using the constraints \eqref{c16}--\eqref{c17},  we obtain
 \berr
  \re^{c_1}\ge \frac{2b_1}{N\lambda\ito\re^{u_1^0+w_1}\ud x},\quad
    \re^{c_2}\ge \frac{2(N-1)b_2}{N\lambda\ito\re^{u_2^0+w_2}\ud x},
 \eerr
 which gives
  \ber
   c_1&\ge& \ln\frac{2b_1}{N}-\ln\lambda-\ln\ito\re^{u_1^0+w_1}\ud x,\label{c35}\\
    c_2&\ge& \ln\frac{2(N-1)b_2}{N}-\ln\lambda-\ln\ito\re^{u_2^0+w_2}\ud x. \label{c36}
  \eer

For any $s\in(0, 1)$, in view of Lemma \ref{lem3} and the Trudinger--Moser inequality \eqref{c32}, we have
 \ber
  \ln\ito\re^{u_1^0+w_1}\ud x&\le& \frac{1-s}{s}\left(\ln\lambda+\ln\frac{N^2}{4[N-1+\kappa]b_1}\right)+\frac1s\ln\ito\re^{su_1^0+sw_1}\ud x\nm\\
   &\le&\frac{s}{16\pi}\|\nabla w_1\|_2^2+\frac1s\ln C_1+\max\limits_{\Omega} u_1^0\nm\\
   &&+ \frac{1-s}{s}\left(\ln\lambda+\ln\frac{N^2}{4[N-1+\kappa]b_1}\right).\label{c37}
 \eer

Similarly, we obtain
\ber
 \ln\ito\re^{u_2^0+w_2}\ud x  &\le&\frac{s}{16\pi}\|\nabla w_2\|_2^2+\frac1s\ln C_1 +\max\limits_{\Omega} u_2^0\nm\\
 &&+ \frac{1-s}{s}\left(\ln\lambda+\ln\frac{N^2}{4b_2[N-1](1+[N-1]\kappa)}\right).\label{c38}
 \eer

Plugging \eqref{c37} and  \eqref{c38} into \eqref{c35}and  \eqref{c36},  respectively, gives
 \ber
    c_1&\ge&-\frac{s}{16\pi}\|\nabla w_1\|_2^2+\ln\frac{N}{2(N-1+\kappa)}-\max\limits_{\Omega}u_1^0\nm\\
    &&-\frac1s\left(\ln\lambda+\ln C_1+\ln\frac{N^2}{4(N-1+\kappa)b_1}\right),\label{c39}\\
    c_2&\ge&-\frac{s}{16\pi}\|\nabla w_2\|_2^2+\ln\frac{N}{2(1+[N-1]\kappa)}-\max\limits_{\Omega}u_2^0\nm\\
    &&-\frac1s\left(\ln\lambda+\ln C_1+\ln\frac{N^2}{4(N-1)(1+[N-1]\kappa)b_2}\right). \label{c40}
 \eer

Hence, inserting \eqref{c39}--\eqref{c40} into \eqref{c33}, we obtain
 \ber
   J(w_1, w_2)&\ge&\left(\frac{\alpha_0(\kappa)}{2}-\frac{b_1s}{16\pi}\right)\|\nabla w_1\|_2^2+\left(\frac{\alpha_0(\kappa)}{2}-\frac{b_2s}{16\pi}\right)\|\nabla w_2\|_2^2\nm\\
     &&-\frac{b_1+b_2}{s}\ln\lambda+b_1\ln\frac{N}{2(N-1+\kappa)}+b_2\ln\frac{N}{2(1+[N-1]\kappa)}\nm\\
     &&-b_1\max\limits_{\Omega}u_1^0-b_2\max\limits_{\Omega}u_2^0-\frac{b_1+b_2}{s}\ln C_1\nm\\
     &&-\frac1s\left( b_1\ln\frac{N^2}{4(N-1+\kappa)b_1}+b_2\ln\frac{N^2}{4(N-1)(1+[N-1]\kappa)b_2}\right). \label{c41}
 \eer

Thus, we get \eqref{c32'} by taking $s>0$ sufficiently small in \eqref{c41}. Then the proof of Lemma \ref{lem4} is complete.

Noting that the functional $J$ is weakly lower semi-continuous  in $\mathcal{A}$, then, using  Lemma \ref{lem4}, we infer that
 $J$ admits a minimizer in $\mathcal{A}$.  In the sequel,  we show that the minimizer  of $J$ lies in the interior of $\mathcal{A}$ when $\lambda$
 is sufficiently large.

 \begin{lemma}\label{lem5}
  On the boundary of $\mathcal{A}$,  there exists a constant $C_4>0$ independent of $\lambda$ such that
   \ber
    \inf\limits_{(w_1, w_2)\in\partial\mathcal{A}} J(w_1, w_2) &\ge& \frac{N|\Omega|}{2(N-1)}\lambda-C_4(\ln\lambda+\sqrt{\lambda}+1).\label{c42}
   \eer
 \end{lemma}

{\bf Proof.} \quad  By the definition of $\mathcal{A}$, we see that  on the boundary of $\mathcal{A}$ there hold
 \ber
   \left(\ito\re^{u_1^0+w_1}\ud x\right)^2&=& \frac{4b_1(N-1+\kappa)}{N^2\lambda}\ito\re^{2u_1^0+2w_1}\ud x,\label{c43}
   \eer
   or
   \ber
   \left(\ito\re^{u_2^0+w_2}\ud x\right)^2&=&  \frac{4b_2(N-1)(1+[N-1]\kappa)}{N^2\lambda}\ito\re^{2u_2^0+2w_2}\ud x.\label{c44}
 \eer

If \eqref{c43} holds,  in view of \eqref{c28} and the H\"{o}lder inequality, we have
\ber
  \re^{c_1}\ito\re^{u_1^0+w_1}\ud x&\le& \frac{\left(1+[N-1]\kappa\right)\left(\ito\re^{u_1^0+w_1}\ud x\right)^2}{N\kappa\ito\re^{2u_1^0+2w_1}\ud x}\nm\\
  &&+\frac{(\kappa-1)\ito\re^{u_1^0+w_1}\ud x\ito\re^{u_2^0+w_2}\ud x\ito\re^{u_1^0+u_2^0+w_1+w_2}\ud x}{N\kappa\ito\re^{2u_1^0+2w_1}\ud x\ito\re^{2u_2^0+2w_2}\ud x}\nm\\
  &\le& \frac{4(N-1+\kappa)(1+[N-1]\kappa)b_1}{N^3\kappa\lambda}+\frac{2(\kappa-1)}{N^2\kappa}\sqrt{\frac{(N-1+\kappa)b_1|\Omega|}{\lambda}}.\label{c45}
\eer

Hence, from Lemma \ref{lem2} and \eqref{c45}, we have
  \ber
   \frac\lambda2\left(N\ito\left[1-\re^{u_1^0+v_1}\right]\ud x+\frac{N}{N-1}\ito\left[1-\re^{u_2^0+v_2}\right]\ud x\right)\ge \frac{N|\Omega|}{2}\lambda-C_5(\sqrt{\lambda}+1), \label{c46}
  \eer
  where $C_5$ is a positive constant independent of $\lambda$.

Similarly, if \eqref{c44} holds,  we can conclude that
 \ber
  \frac\lambda2\left(N\ito\left[1-\re^{u_1^0+v_1}\right]\ud x+\frac{N}{N-1}\ito\left[1-\re^{u_2^0+v_2}\right]\ud x\right)\ge \frac{N|\Omega|}{2(N-1)}\lambda-C_6(\sqrt{\lambda}+1), \label{c47}
  \eer
 where $C_6$ is a positive constant independent of $\lambda$.

Now,  using \eqref{c46}--\eqref{c47}, and   estimating    $c_1, c_2$ as that  in Lemma \ref{lem4}, we can obtain \eqref{c42}.      Then  Lemma \ref{lem5}  follows.

To proceed further, we use the approach of \cite{taran96} need to find the test functions, which lie in the interior of $\mathcal{A}$.

It was shown in \cite{taran96} that, for $\mu>0$ sufficiently large, the problem
 \berr
  \Delta v=\mu\re^{u_i^0+v}(\re^{u_i^0+v}-1)+\frac{4\pi n_i}{|\Omega|} \quad  \text{in}\quad \Omega
 \eerr
admit solutions $v_i^\mu  , \, i=1,2$ such that
  $u_i^0+v_i^\mu<0$ in $\Omega$, $c_i^\mu=\frac{}{}\ito v_i^\mu\ud x \to 0$ and   $w_i^\mu=v_i^\mu-c_i^\mu\to -u_i^0$ pointwise and a.e  as $\mu\to +\infty$.

 Since $\re^{u_i^0+w_i^\mu}\in L^\infty(\Omega)$, $i=1,2$ we have
 \[ \re^{u_i^0+w_i^\mu}\to 1\quad \text{strongly in }\quad  L^p(\Omega)\quad \text{for any } \quad p\ge1 \]
as $\mu\to +\infty$.  In particular, we have
 \berr
  \ito \re^{2u_i^0+2w_i^\mu}\ud x\to |\Omega|, \quad i=1, 2
 \eerr
 as $\mu\to +\infty$.

Thus, for $\lambda_0$ large and for fixed $\vep\in(0, 1)$, we may find $\mu_\vep\gg1$, such that $(w_1^{\mu_\vep}, w_2^{\mu_\vep})\in int\mathcal{A}$ for
every $\lambda>\lambda_0$, and there holds
 \ber
 &&\frac{N^2\kappa|\Omega|^2}{(N-1+\kappa)\big(1+[N-1]\kappa\big)\ito\re^{2u_1^0+2w_1^{\mu_\vep}}\ud x\ito\re^{2u_2^0+2w_2^{\mu_\vep}}\ud x-(N-1)(\kappa-1)^2|\Omega|^2}\nm\\
 &&\ge 1-\vep. \label{c48}
 \eer

Using  Jensen's inequality,  \eqref{c25'},  and \eqref{c19}--\eqref{c20}, we obtain
 \ber
  \re^{c_1(w_1^{\mu_\vep},  w_2^{\mu_\vep})}&\ge&\frac{P_1\big(w_1^{\mu_\vep}, w_2^{\mu_\vep}, \re^{c_2(w_1^{\mu_\vep},  w_2^{\mu_\vep})}\big)}{2(N-1+\kappa)\ito\re^{2u_1^0+2w_1^{\mu_\vep}}\ud x}\nm\\&&
  \times\left(1+\sqrt{1-\frac{4(N-1+\kappa)b_1\ito\re^{2u_1^0+2w_1}\ud x}{\lambda P_1^2\big(w_1^{\mu_\vep}, w_2^{\mu_\vep}, \re^{c_2(w_1^{\mu_\vep},  w_2^{\mu_\vep})}\big)}}\right)\nm\\
  &\ge&\frac{P_1\big(w_1^{\mu_\vep}, w_2^{\mu_\vep}, \re^{c_2(w_1^{\mu_\vep},  w_2^{\mu_\vep})}\big)}{(N-1+\kappa)\ito\re^{2u_1^0+2w_1^{\mu_\vep}}\ud x}-\frac{2b_1}{\lambda P_1\big(w_1^{\mu_\vep}, w_2^{\mu_\vep}, \re^{c_2(w_1^{\mu_\vep},  w_2^{\mu_\vep})}\big)}\nm\\
  &\ge& \frac{\big(N+[\kappa-1]\re^{c_2(w_1^{\mu_\vep},  w_2^{\mu_\vep})}\big)|\Omega|}{(N-1+\kappa)\ito\re^{2u_1^0+2w_1^{\mu_\vep}}\ud x}-\frac{2b_1}{N\lambda|\Omega|}.\label{c49}
 \eer

Analogously, we have
 \ber
   \re^{c_2(w_1^{\mu_\vep},  w_2^{\mu_\vep})}
  &\ge& \frac{\big(N+[N-1][\kappa-1]\re^{c_1(w_1^{\mu_\vep},  w_2^{\mu_\vep})}\big)|\Omega|}{(1+[N-1]\kappa)\ito\re^{2u_2^0+2w_2^{\mu_\vep}}\ud x}-\frac{2(N-1)b_2}{N\lambda|\Omega|}.\label{c50}
 \eer

Plugging \eqref{c50} into \eqref{c49}, we see that
 \berr
  &&\re^{c_1(w_1^{\mu_\vep},  w_2^{\mu_\vep})}\nm\\
  &&\ge \frac{N|\Omega|}{(N-1+\kappa)\ito\re^{2u_1^0+2w_1^{\mu_\vep}}\ud x}-\frac{2b_1}{N\lambda|\Omega|}\nm\\
  &&+ \frac{(\kappa-1)|\Omega|}{(N-1+\kappa)\ito\re^{2u_1^0+2w_1^{\mu_\vep}}\ud x}\left(\frac{\big(N+[N-1][\kappa-1]\re^{c_1(w_1^{\mu_\vep},  w_2^{\mu_\vep})}\big)|\Omega|}{(1+[N-1]\kappa)\ito\re^{2u_2^0+2w_2^{\mu_\vep}}\ud x}-\frac{2(N-1)b_2}{N\lambda|\Omega|}\right)\nm\\
  &&\ge \frac{N|\Omega|\left((1+[N-1]\kappa)\ito\re^{2u_2^0+2w_2^{\mu_\vep}}\ud x+(\kappa-1)|\Omega|\right)+(N-1)(\kappa-1)^2|\Omega|^2\re^{c_1(w_1^{\mu_\vep},  w_2^{\mu_\vep})}}{(N-1+\kappa)(1+[N-1]\kappa)\ito\re^{2u_1^0+2w_1^{\mu_\vep}}\ud x\ito\re^{2u_2^0+2w_2^{\mu_\vep}}\ud x}\nm\\
  &&-\frac{2}{N\lambda|\Omega|}\left(b_1+\frac{[N-1][\kappa-1]b_2}{N-1+\kappa}\right)\nm\\
  &&\ge \frac{N^2\kappa|\Omega|^2+(N-1)(\kappa-1)^2|\Omega|^2\re^{c_1(w_1^{\mu_\vep},  w_2^{\mu_\vep})}}{(N-1+\kappa)(1+[N-1]\kappa)\ito\re^{2u_1^0+2w_1^{\mu_\vep}}\ud x\ito\re^{2u_2^0+2w_2^{\mu_\vep}}\ud x}\nm\\
  &&-\frac{2}{N\lambda|\Omega|}\left(b_1+\frac{[N-1][\kappa-1]b_2}{N-1+\kappa}\right),
 \eerr
which implies
 \ber
  &&\re^{c_1(w_1^{\mu_\vep},  w_2^{\mu_\vep})}\nm\\
  &&\ge\frac{N^2\kappa|\Omega|^2}{(N-1+\kappa)(1+[N-1]\kappa)\ito\re^{2u_1^0+2w_1^{\mu_\vep}}\ud x\ito\re^{2u_2^0+2w_2^{\mu_\vep}}\ud x-(N-1)(\kappa-1)^2|\Omega|^2}\nm\\
  &&-\frac{(N-1+\kappa)(1+[N-1]\kappa)\ito\re^{2u_1^0+2w_1^{\mu_\vep}}\ud x\ito\re^{2u_2^0+2w_2^{\mu_\vep}}\ud x}{(N-1+\kappa)(1+[N-1]\kappa)\ito\re^{2u_1^0+2w_1^{\mu_\vep}}\ud x\ito\re^{2u_2^0+2w_2^{\mu_\vep}}\ud x-(N-1)(\kappa-1)^2|\Omega|^2}\nm\\
  &&\times\frac{2}{N\lambda|\Omega|}\left(b_1+\frac{[N-1][\kappa-1]b_2}{N-1+\kappa}\right)\nm\\
  &&\ge\frac{N^2\kappa|\Omega|^2}{(N-1+\kappa)(1+[N-1]\kappa)\ito\re^{2u_1^0+2w_1^{\mu_\vep}}\ud x\ito\re^{2u_2^0+2w_2^{\mu_\vep}}\ud x-(N-1)(\kappa-1)^2|\Omega|^2}\nm\\
  &&-\frac{2(1+[N-1]\kappa)\left([N-1+\kappa]b_1+[N-1][\kappa-1]b_2\right)}{N^3\lambda\kappa|\Omega|}.\label{c51}
 \eer
Similarly, we have
   \ber
  &&\re^{c_2(w_1^{\mu_\vep},  w_2^{\mu_\vep})}\nm\\
  &&\ge\frac{N^2\kappa|\Omega|^2}{(N-1+\kappa)(1+[N-1]\kappa)\ito\re^{2u_1^0+2w_1^{\mu_\vep}}\ud x\ito\re^{2u_2^0+2w_2^{\mu_\vep}}\ud x-(N-1)(\kappa-1)^2|\Omega|^2}\nm\\
  &&-\frac{2(N-1+\kappa)\left\{(1+[N-1]\kappa)b_1+[N-1][\kappa-1]b_2\right\}}{N^3\lambda\kappa|\Omega|}.\label{c52}
 \eer

Hence, from \eqref{c47}, \eqref{c51} and \eqref{c52}, we conclude that, for all $\lambda>\lambda_0$,
 \berr
   \re^{c_1(w_1^{\mu_\vep},  w_2^{\mu_\vep})}&\ge&1-\vep-\frac{2(1+[N-1]\kappa)\left([N-1+\kappa]b_1+[N-1][\kappa-1]b_2\right)}{N^3\lambda\kappa|\Omega|},\\
   \re^{c_2(w_1^{\mu_\vep},  w_2^{\mu_\vep})}&\ge&1-\vep-\frac{2(N-1+\kappa)\left\{(1+[N-1]\kappa)b_1+[N-1][\kappa-1]b_2\right\}}{N^3\lambda\kappa|\Omega|}.
 \eerr

Consequently, we have
 \ber
 &&\ito\left(1-\re^{c_1(w_1^{\mu_\vep},  w_2^{\mu_\vep})}\re^{u_1^0+w_1^{\mu_\vep}}\right)\nm\\
 &&\le|\Omega|\vep-\frac{2(1+[N-1]\kappa)\left([N-1+\kappa]b_1+[N-1][\kappa-1]b_2\right)}{N^3\lambda\kappa},\label{c53}\\
 &&\ito\left(1-\re^{c_2(w_1^{\mu_\vep},  w_2^{\mu_\vep})}\re^{u_2^0+w_2^{\mu_\vep}}\right)\nm\\
 &&\le|\Omega|\vep-\frac{2(N-1+\kappa)\left\{(1+[N-1]\kappa)b_1+[N-1][\kappa-1]b_2\right\}}{N^3\lambda\kappa},\label{c54}
 \eer
 for all $\lambda>\lambda_0$.

\begin{lemma}\label{lem6}
  As $\lambda>0$ is  sufficiently large, there holds
   \ber
    J(w_1^{\mu_\vep},  w_2^{\mu_\vep})-\inf\limits_{(w_1, w_2)\in\partial\mathcal{A}}J(w_1, w_2)<-1.\label{c55}
   \eer
\end{lemma}

{\bf Proof.} \quad  Using \eqref{c25'}, \eqref{c53} and \eqref{c54}, we infer that, for any small $\vep>0$, there exists a positive constant $C_\vep$ such
that
\ber
 J(w_1^{\mu_\vep},  w_2^{\mu_\vep})\le \frac{N^2\lambda|\Omega|\vep}{2(N-1)}+C_\vep. \label{c56}
\eer
Thus, in view of Lemma \ref{lem5}, we have
  \ber
   J(w_1^{\mu_\vep},  w_2^{\mu_\vep})-\inf\limits_{(w_1, w_2)\in\partial\mathcal{A}} J(w_1, w_2)\le  \frac{N|\Omega|\lambda}{2(N-1)}(N\vep-1)+C(\ln\lambda+\sqrt{\lambda}+1),\label{c57}
  \eer
where $C$ is a positive constant independent of $\lambda$.

Then, taking $\vep=\frac{1}{2N}$, and $\lambda$ sufficiently large in \eqref{c57}, we conclude  \eqref{c55}.

Now    from Lemma \ref{lem4} and Lemma \ref{lem6} we infer the following corollary.
 \begin{corollary}
  There exists $\tilde{\lambda}>0$ such that, for every $\lambda>\tilde{\lambda}$, the functional $J$ achieves its minimum at a point
  $(w_{1, \lambda}, w_{2, \lambda})$, which belongs to the interior of $\mathcal{A}$. Moreover,  $(v_{1, \lambda}, v_{2, \lambda})$,  defined by
   \ber
    v_{i, \lambda}=w_{i, \lambda}+c_i(w_{1, \lambda}, w_{2, \lambda}), \quad i=1, 2, \label{c58}
   \eer
 \end{corollary}
is a critical point of the functional $I$ in $W^{1, 2}(\Omega)\times W^{1, 2}(\Omega)$, namely, a weak solution of \eqref{c1}--\eqref{c2}.

Next we study the behavior of the solution given above.

\begin{lemma}\label{lem7}
   Let $(v_{1, \lambda}, v_{2, \lambda})$ be the solution of \eqref{c1}--\eqref{c2} given by \eqref{c58}. There holds
    \ber
     \re^{u_i^0+v_{i, \lambda}}\to 1 \quad \text{as} \quad \lambda\to +\infty, \quad i=1, 2,\label{c59}
   \eer
    pointwise a.e. in $\Omega$ and in $L^p(\Omega)$ for any $p\ge1$.
   Moreover,  $(v_{1, \lambda}, v_{2, \lambda})$  is a local minimizer of the functional  $I$ in  $W^{1, 2}(\Omega)\times W^{1, 2}(\Omega)$.
\end{lemma}

{\bf Proof.} \quad Using   \eqref{c8} and   similar estimates as in Lemma \ref{lem4}, for any $\lambda>\tilde{\lambda}$, we  infer that  there
exists a positive constant $C$  independent of $\lambda$ such that
  \ber
   J(w_{1, \lambda}, w_{2, \lambda})&\ge&  \frac{\alpha_0(\kappa^{-1})\lambda}{2}\left\{\ito\left(\re^{u_1^0+v_{1, \lambda}}-1\right)^2\ud x+\ito\left(\re^{u_2^0+v_{2, \lambda}}-1\right)^2\ud x\right\}\nm\\
   &&-C(\ln\lambda+1),\label{c60}
   \eer
   where $\alpha_0(\kappa^{-1})$ is a positive constant defined by \eqref{a140}.
Hence, it follows from \eqref{c60} and  \eqref{c56} that
  \ber
   \ito\left(\re^{u_i^0+v_{i, \lambda}}-1\right)^2\ud x\to 0, \quad \text{as} \quad \lambda\to +\infty, \quad i=1, 2.\label{c61}
  \eer
  In view of Lemma \ref{lem1}, we have $\re^{u_i^0+v_{i, \lambda}}<1, \, i=1,2$. Then, we conclude  \eqref{c59} by the dominated convergence theorem.

Next, we show that  $(v_{1, \lambda}, v_{2, \lambda})$  is a local minimizer of the functional  $I$ in  $W^{1, 2}(\Omega)\times W^{1, 2}(\Omega)$.

  By a direct computation, for any $(w_1, w_2)\in\mathcal{A}$ and the corresponding $(c_1, c_2)$ given by \eqref{c19}--\eqref{c20}, we obtain
  \berr
    \partial_{c_1}I(w_1+c_1(w_1, w_2), w_2+c_2(w_1,w_2)) =0= \partial_{c_2}I(w_1+c_1(w_1, w_2), w_2+c_2(w_1,w_2))
   \eerr
and
 \ber
    &&\partial^2_{c^2_1}I(w_1+c_1(w_1, w_2), w_2+c_2(w_1,w_2))\nm\\
    &&= \lambda\left(2[N-1+\kappa]\re^{2c_1}\ito\re^{2u_1^0+2w_1}\ud x-\re^{c_1}P_1(w_1, w_2, \re^{c_2})\right)\nm\\
    &&=\lambda\left\{\left(N\ito\re^{u_1^0+v_1}\ud x+[\kappa-1]\ito\re^{u_1^0+u_2^0+v_1+v_2}\ud x\right)^2\right.\nm\\
    &&\quad \left.-\frac{4(N-1+\kappa)b_1}{\lambda}\ito\re^{2u_1^0+2v_1}\ud x\right\}^{\frac12},\label{b59}
    \eer
    \ber
    &&\partial^2_{c^2_2}I(w_1+c_1(w_1, w_2), w_2+c_2(w_1,w_2)) \nm\\
    &&= \lambda\left(2\left[\frac{1}{N-1}+\kappa\right]\re^{2c_2}\ito\re^{2u_2^0+2w_2}\ud x-\re^{c_2}P_2(w_1, w_2, \re^{c_1})\right)\nm\\
    &&=\lambda\left\{\left(\frac{N}{N-1}\ito\re^{u_2^0+v_2}\ud x+[\kappa-1]\ito\re^{u_1^0+u_2^0+v_1+v_2}\ud x\right)^2\right.\nm\\
    &&\quad \left.-\frac{4(1+[N-1]\kappa)b_2}{(N-1)\lambda}\ito\re^{2u_2^0+2v_2}\ud x\right\}^{\frac12},\label{b60}\\
    &&\partial^2_{c_1c_2}I(w_1+c_1(w_1, w_2), w_2+c_2(w_1,w_2)) \nm\\
    &&= \lambda(1-\kappa)\re^{c_1}\re^{c_2} \ito\re^{u_1^0+u_2^0+w_1+w_2}\ud x\nm\\
    &&=\lambda(1-\kappa)\ito\re^{u_1^0+u_2^0+v_1+v_2}\ud x.\label{b61}
\eer

 If  $(w_1, w_2)$ belongs to  the interior of $\mathcal{A}$,  then we can use strict inequalities in the constraints \eqref{c16}--\eqref{c17} to
 get
  \berr
   \partial^2_{c^2_1}I(w_1+c_1(w_1, w_2), w_2+c_2(w_1,w_2))>\lambda{(\kappa-1)}\ito\re^{u_1^0+u_2^0+v_1+v_2}\ud  x,\\
   \partial^2_{c^2_2}I(w_1+c_1(w_1, w_2), w_2+c_2(w_1,w_2))>\lambda{(\kappa-1)}\ito\re^{u_1^0+u_2^0+v_1+v_2}\ud  x.
  \eerr
Thus,  we conclude  that, if $(w_1, w_2)$ is an interior point of  $ \mathcal{A}$ then the Hessian matrix of $I(w_1+c_1, w_2+c_2)$ with respect to $(c_1, c_2)$
is strictly positive definite at $(c_1(w_1, w_2), c_2(w_1, w_2))$. We apply such property,
near the critical point $(v_{1, \lambda}, v_{2, \lambda})$. Indeed, by continuity, for $\delta>0$ sufficiently small, we can ensure that, if
 $(v_1, v_2)=(w_1+c_1, w_2+c_2)$  satisfies:
  \[ \|v_1-v_{1, \lambda}\|+\|v_2-v_{2, \lambda}\|\le \delta, \]
  then  $(w_1, w_2)$ belongs to the interior of $\mathcal{A}$  and
  \berr
   I(v_1, v_2)= I(w_1+c_1, w_2+c_2)&\ge& I(w_{1, \lambda}+c_1(w_{1,\lambda}, w_{2, \lambda}), w_{2,\lambda}+c_2(w_{1, \lambda}, w_{2, \lambda}))\\
   &=&J(w_{1,\lambda}, w_{2, \lambda})= I(v_{1, \lambda}, v_{2, \lambda}).
  \eerr
Hence, $(v_{1, \lambda}, v_{2, \lambda})$  is  a  local minimizer for $ I$ in  $W^{1, 2}(\Omega)\times W^{1, 2}(\Omega)$.
Then the proof of Lemma \ref{lem7} is complete.

\subsection{A second  solution}
 In this subsection, via mountain-pass theorem, we find a second critical point of the functional $I$, which gives a second solution of \eqref{c1}--\eqref{c2}.

For  this purpose, we  show that the functional $I$ satisfies the P-S condition.
\begin{lemma} \label{lem8}
Every sequence $(v_{1, n}, v_{2, n})\in W^{1, 2}(\Omega)\times W^{1, 2}(\Omega)$ satisfies
\ber
 I(v_{1, n}, v_{2, n})\to a_0\quad \text{as} \quad n\to +\infty,\label{c62} \\
 \|I'(v_{1, n}, v_{2, n})\|_*\to 0 \quad \text{as} \quad n\to +\infty,\label{c63}
\eer
admits a strongly convergent subsequence in $W^{1, 2}(\Omega)\times W^{1, 2}(\Omega)$, where  $a_0$ is a constant and $\|\cdot\|_*$ denotes the
norm of the dual space of  $W^{1, 2}(\Omega)\times W^{1, 2}(\Omega)$.
\end{lemma}

{\bf Proof.}\quad  Denote   $\vep_n=\|I'(v_{1, n}, v_{2, n})\|_*$,  we have  $\vep_n\to 0$ as $n\to +\infty$. For any $(\psi_1, \psi_2)\in W^{1, 2}(\Omega)\times W^{1, 2}(\Omega)$,
we obtain
\ber
 &&(I'(v_{1, n}, v_{2, n}))(\psi_1, \psi_2)\nm\\
 &&= \left(N-1+\frac1\kappa\right)\ito\nabla v_{1, n}\cdot\nabla\psi_1\ud x+\left(\frac{1}{N-1}+\frac1\kappa\right)\ito \nabla v_{2, n}\nabla \psi_2\ud x
 \nm\\&&+\left(1-\frac1\kappa\right)\ito(\nabla v_{2, n}\cdot\nabla\psi_1+\nabla v_{1, n}\cdot\nabla\psi_2)\ud x\nm\\
 &&+\lambda\ito\left(\left[N-1+\kappa\right]\re^{u_1^0+v_{1, n}}\left[\re^{u_1^0+v_{1, n}}-1\right]+(1-\kappa)\re^{u_1^0+v_{1, n}}\left[\re^{u_2^0+v_{2, n}}-1\right]\right)\psi_1\ud x\nm\\
 &&+\lambda\ito\left(\left[\frac{1}{N-1}+\kappa\right]\re^{u_2^0+v_{2, n}}\left[\re^{u_2^0+v_{2, n}}-1\right]+(1-\kappa)\re^{u_2^0+v_{2, n}}\left[\re^{u_1^0+v_{1, n}}-1\right]\right)\psi_2\ud x\nm\\
 &&+\frac{b_1}{|\Omega|}\ito \psi_1\ud x+\frac{b_2}{|\Omega|}\ito \psi_2\ud x\label{c64}
\eer
and
\be
|(I'(v_{1, n}, v_{2, n}))(\psi_1, \psi_2)|\le \vep_n(\|\psi_1\|+\|\psi_2\|)\label{c64'}
\ee

 Taking $(\psi_1, \psi_2)=(1, 1)$ in \eqref{c64}, we  find
\ber
 &&(I'(v_{1, n}, v_{2, n}))(1, 1)\nm\\
 &&=\lambda\ito\left(\left[N-1+\kappa\right]\re^{u_1^0+v_{1, n}}\left[\re^{u_1^0+v_{1, n}}-1\right]+(1-\kappa)\re^{u_1^0+v_{1, n}}\left[\re^{u_2^0+v_{2, n}}-1\right]\right)\ud x\nm\\
 &&+\lambda\ito\left(\left[\frac{1}{N-1}+\kappa\right]\re^{u_2^0+v_{2, n}}\left[\re^{u_2^0+v_{2, n}}-1\right]+(1-\kappa)\re^{u_2^0+v_{2, n}}\left[\re^{u_1^0+v_{1, n}}-1\right]\right)\ud x\nm\\
 &&+b_1+b_2\nm\\
 &&=\lambda\left\{\ito\left(\left[N-1+\kappa\right]\left[\re^{u_1^0+v_{1, n}}-1\right]^2+\left[\frac{1}{N-1}+\kappa\right]\left[\re^{u_2^0+v_{2, n}}-1\right]^2\right)\ud x\right.\nm\\
 &&+ 2(1-\kappa)\ito\left(\re^{u_1^0+v_{1, n}}-1\right)\left(\re^{u_2^0+v_{2, n}}-1\right)\ud x\nm\\
 &&\left.+\ito\left(N\left[\re^{u_1^0+v_{1, n}}-1\right]+\frac{N}{N-1}\left[\re^{u_2^0+v_{2, n}}-1\right]\right)\ud x\right\}+b_1+b_2.\label{c65}
\eer
 Noting that the matrix $A(N, \kappa^{-1})$ defined by \eqref{11} is positive definite,  then from \eqref{c65} and \eqref{c64'}  we infer that
 \berr
\ito\left(\re^{u_1^0+v_{1, n}}-1\right)^2\ud x+\ito\left(\re^{u_2^0+v_{2, n}}-1\right)^2\ud x\le C
 \eerr
for some positive constant $C$, which implies
 \be
   \ito\re^{2u_1^0+2v_{1, n}}\ud x+\ito\re^{2u_2^0+2v_{2, n}}\ud x\le C. \label{c66}
 \ee
 Here and what follows we use $C$ to denote a generic positive constant independent of $n$.
Since $v_{i, n}\in W^{1, 2}(\Omega)$,  we have the following decomposition
\berr
 v_{i, n}=w_{i, n}+c_{i, n}, \quad w_{i, n}\in \dot{W}^{1, 2}(\Omega), \quad c_{i, n}\in \mathbb{R}.
\eerr
From \eqref{c66} we conclude that $c_{i, n}$  is bounded from above.

 Noting the matrix $A(N, \kappa)$ and $A(N, \kappa^{-1})$ defined by   \eqref{11} are both positive definite,  we estimate  $I(v_{1, n}, v_{2, n})$ as
\ber
 &&I(v_{1, n}, v_{2, n})\ge \frac{\alpha_0(\kappa)}{2}\left(\|\nabla w_{1, n}\|_2^2+\|\nabla w_{2, n}\|_2^2\right)\nm\\
 &&+\frac{\alpha_0(\kappa^{-1})\lambda}{2}\left(\ito\left[\re^{u_1^0+v_{1, n}}-1\right]^2\ud x\ito\left[\re^{u_2^0+v_{2, n}}-1\right]^2\ud x\right)+ b_1c_{1, n}+ b_2c_{2, n},\label{c67}
\eer
where $\alpha_0(\kappa)$   and  $\alpha_0(\kappa^{-1})$  are a positive constants defined by \eqref{a140}.

Let $\varphi_n\equiv Nw_{1,n}+\frac{N}{N-1}w_{2, n}$  and $\varphi_n^+\equiv\max\{\varphi_n, 0\}$. Then, taking  $ (\psi_1, \psi_2)=(\varphi_n^+, \varphi_n^+)$ in \eqref{c64}, we get
\ber
 &&\|\nabla \varphi_n^+\|_2^2+\lambda\ito\left(\sqrt{N-1+\kappa}\re^{u_1^0+v_{1,n}}-\sqrt{\frac{1}{N-1}+\kappa}\re^{u_2^0+v_{2, n}}\right)^2\varphi_n^+\ud x\nm\\
 &&+2\lambda\left(\sqrt{(\kappa-1)^2+\frac{N^2\kappa}{N-1}}-[\kappa-1]\right)\ito \re^{u_1^0+u_2^0+v_{1, n}+v_{2, n}}\varphi_n^+\ud x\nm\\
 &&\le C(\|\varphi_n^+\|_2+\vep_n\|\varphi_n^+\|).\label{c68}
\eer

Therefore, using the Poincar\'{e} inequality in  \eqref{c68}, we find that
  \ber
  \ito \re^{u_1^0+u_2^0+v_{1, n}+v_{2, n}}\varphi_n^+\ud x\le C(\|\nabla w_{1, n}\|_2+\|\nabla w_{2, n}\|_2). \label{c69}
  \eer

Let $(\psi_1, \psi_2)=(w_{1, n}, w_{2, n})$ in \eqref{c64},  and noting that the matrix $A(N, \kappa)$  defined by \eqref{11} is positive definite,  we have
\ber
 &&(I'(v_{1, n}, v_{2, n}))(w_{1, n}, w_{2, n})\nm\\
 &&=\left(N-1+\frac1\kappa\right)\|\nabla w_{1, n}\|_2^2+\left(\frac{1}{N-1}+\frac1\kappa\right)\left\|\nabla w_{2, n}\right\|_2^2+2\left(1-\frac1\kappa\right)\ito\nabla w_{1, n}\cdot \nabla w_{2, n}\ud x\nm\\
 &&+\lambda\left\{\left[N-1+\kappa\right]\ito\re^{2u_1^0+2v_{1, n}}w_{1, n}\ud x + \left[\frac{1}{N-1}+\kappa\right]\ito\re^{2u_2^0+v_{2, n}}w_{2, n}\ud x\right.\nm\\
 &&\left.+(1-\kappa)\ito\re^{u_1^0+2v_{1, n}+u_2^0+v_{2, n}}(w_{1, n}+w_{2, n})\ud x\right.\nm\\
 &&\left.-N\ito\re^{u_1^0+v_{1, n}}w_{1, n}\ud x-\frac{N}{N-1}\ito\re^{2u_2^0+2v_{2, n}}w_{2, n}\ud x\right\}\nm\\
 &&\ge \alpha_0(\kappa)\left(\|\nabla w_{1, n}\|_2^2+\|\nabla w_{2, n}\|_2^2\right)\nm\\
 &&+\lambda\left\{\left[N-1+\kappa\right]\ito\re^{2u_1^0+2v_{1, n}}w_{1, n}\ud x + \left[\frac{1}{N-1}+\kappa\right]\ito\re^{2u_2^0+v_{2, n}}w_{2, n}\ud x\right.\nm\\
 &&\left.+(1-\kappa)\ito\re^{u_1^0+2v_{1, n}+u_2^0+v_{2, n}}(w_{1, n}+w_{2, n})\ud x\right.\nm\\
 &&\left.-N\ito\re^{u_1^0+v_{1, n}}w_{1, n}\ud x-\frac{N}{N-1}\ito\re^{2u_2^0+2v_{2, n}}w_{2, n}\ud x\right\}, \label{c70}
\eer
where $\alpha_0(\kappa)$  is a  positive constant  defined by \eqref{a140}.
In view of \eqref{c66}, we see that
 \ber
  \left|\ito\re^{u_i^0+v_{i, n}}w_{i, n}\ud x\right|\le C\|w_{i, n}\|_2, \quad i=1, 2. \label{c71}
 \eer

Since we have shown that $c_{i, n}$ is bounded from above,   there holds
 \ber
  \ito \re^{2u_i^0+2v_{i, n}}w_{i, n}\ud x&=&\ito\re^{2u_i^0+2c_{i, n}}\left(\re^{2w_{i, n}}-1\right)w_{i, n}\ud x+\ito\re^{2u_i^0+2c_{i, n}}w_{i, n}\ud x\nm\\
  &\ge&-C\|w_{i, n}\|_2.\label{c72}
 \eer

 It easily follows that
   \berr
    &&\ito\re^{u_1^0+v_{1,n}+u_2^0+v_{2,n}}(w_{1,n}+w_{2,n})\ud x\\
    &&\le \ito\re^{u_1^0+v_{1,n}+u_2^0+v_{2,n}}(w_{1,n}+w_{2,n})_+\ud x\\
    && =\int\limits_{\{w_{1,n}\le0\le w_{2,n}\}} \re^{u_1^0+c_{1,n}}\left(\re^{w_{1,n}}-1\right)\re^{u_2^0+v_{2,n}}(w_{1,n}+w_{2,n})_+\ud x\\
    && \quad+\int\limits_{\{w_{1,n}\le0\le w_{2,n}\}}\re^{u_1^0+c_{1,n}}\re^{u_2^0+v_{2,n}}(w_{1,n}+w_{2,n})_+\ud x\\
    &&\quad +\int\limits_{\{w_{2,n}\le0\le w_{1,n}\}} \re^{u_2^0+c_{2,n}}\left(\re^{w_{2,n}}-1\right)\re^{u_1^0+v_{1,n}}(w_{1,n}+w_{2,n})_+\ud x\\
    &&  \quad+\int\limits_{\{w_{2,n}\le0\le w_{1,n}\}}\re^{u_2^0+c_{2,n}}\re^{u_1^0+v_{1,n}}(w_{1,n}+w_{2,n})_+\ud x\\
    &&\quad +\int\limits_{ \{w_{1,n}>0\}\cap\{w_{2,n}>0\}} \re^{u_1^0+v_{1,n}+u_2^0+v_{2,n}}(w_{1,n}+w_{2,n})_+\ud x\\
    &&\le C\left(\|\nabla w_{1,n}\|_2+\|\nabla w_{2,n}\|_2\right)+\ito\re^{u_1^0+v_{1,n}+u_2^0+v_{2,n}}\varphi_n^+\ud x,
   \eerr
which together with \eqref{c69} imply
  \be
   \ito\re^{u_1^0+v_{1,n}+u_2^0+v_{2,n}}(w_{1,n}+w_{2,n})\ud x\le  C\left(\|\nabla w_{1,n}\|_2+\|\nabla w_{2,n}\|_2\right).\label{c73}
  \ee
Now from \eqref{c70}--\eqref{c73}, we see that
\ber
 \|\nabla w_{1,n}\|_2+\|\nabla w_{2,n}\|_2\le C.\label{c74}
\eer

Noting that we have shown that $\{c_{i, n}\}$ is bounded from above, by \eqref{c74}, \eqref{c62} and \eqref{c67}, we infer that $c_{i, n} $ is also bounded from below, $\, i=1, 2$.
Hence, using \eqref{c74} again, we conclude that $\{v_{i, n}\}$ is uniformly bounded in $W^{1, 2}(\Omega)$, $i=1, 2$.

Therefore, up to a subsequence, there exists $v_i\in W^{1, 2}(\Omega)$, such that $v_{i, n}\to v_i$  weakly in $W^{1, 2}(\Omega)$, strongly in
$L^p(\Omega)$  for any $p\ge1$, pointwise a.e. in $\Omega$, and  $\re^{u_i^0+v_{i, n}}\to \re^{u_i^0+v_i}$ in $L^p(\Omega)$  for any $p\ge 1$, as $n\to+\infty$, $i=1, 2$.

Hence we see that  $(v_1, v_2)$ is a critical point for the functional $I$.  From the above convergence results we obtain
\ber
&&\alpha_0(\kappa)\left(\|\nabla (v_{1, n}-v_1)\|_2^2+\left\|\nabla (v_{2, n}-v_2)\right\|_2^2\right)\nm\\
&&\le \left(N-1+\frac1\kappa\right)\|\nabla (v_{1, n}-v_1)\|_2^2+\left(\frac{1}{N-1}+\frac1\kappa\right)\left\|\nabla (v_{2, n}-v_2)\right\|_2^2\nm\\
&&+2\left(1-\frac1\kappa\right)\ito\nabla (v_{1, n}-v_1)\cdot \nabla (v_{2, n}-v_2)\ud x\nm\\
&&=\big(I'(v_{1, n}, v_{2, n})-I'(v_1, v_2)\big)(v_{1, n}-v_1, v_{2, n}-v_2)+o(1)\to 0\quad \text{as}\quad n\to+\infty,\label{c74a}
\eer
where $\alpha_0(\kappa)$  is a  positive constant  defined by \eqref{a140}.

Then we conclude from the estimate \eqref{c74a} that   $(v_{1, n}, v_{2, n})\to (v_1, v_2)$ strongly in $W^{1, 2}(\Omega)\times W^{1, 2}(\Omega)$ as $n\to +\infty$.
Then Lemma \ref{lem8} follows.

To find a second solution of  \eqref{c1}--\eqref{c2}, noting that we have proved  that  $(v_{1, \lambda}, v_{2, \lambda})$  given in \eqref{c58} is a local minimizer of the functional $I$,
we only need to consider the following two cases.

Case 1.   $(v_{1, \lambda}, v_{2, \lambda})$  is a degenerate minimum.   In other words, for any sufficiently small $\delta>0$,
\berr
 \inf\limits_{\|v_1-v_{1, \lambda}\|+\|v_2-v_{2, \lambda}\|=\delta}I(v_1, v_2)=I(v_{1, \lambda}, v_{2, \lambda}).
\eerr
Thus, we conclude from  Corollary 1.6 of \cite{gho}) that there is  a one parameter family of degenerate local minimizer of the functional $I$. Automatically,  a second solution of  \eqref{c1}--\eqref{c2} for  this case can be obtained.

 Case 2.    $(v_{1, \lambda}, v_{2, \lambda})$  is a  strict local minimum.  That is, for any sufficiently small $\delta>0$,
 there holds
  \ber
   I(v_{1, \lambda}, v_{2, \lambda})<\inf\limits_{\|v_1-v_{1, \lambda}\|+\|v_2-v_{2, \lambda}\|=\delta}I(v_1, v_2)\equiv\gamma_0.\label{c75}
  \eer

We observe  that
 \berr
  I(v_{1, \lambda}-\xi, v_{2, \lambda}-\xi)\to -\infty\quad \text{as}\quad \xi\to +\infty.
 \eerr

Hence,  for a sufficiently large $\xi_0>1$, let
\berr
 \tilde{v}_i= v_{i, \lambda}-\xi_0, \quad i=1, 2,
\eerr
we can obtain
 \ber
  \|\tilde{v}_1-v_{1, \lambda}\|+\|\tilde{v}_2-v_{2, \lambda}\|>\delta\label{c76}
 \eer
 and
 \ber
  I(\tilde{v}_1, \tilde{v}_2)<I(v_{1, \lambda}, v_{2, \lambda})-1\label{c77}
 \eer

Now we  introduce the paths
\berr
\mathcal{P}=\Big\{ \Gamma(t)\Big| \Gamma\in C\left([0, 1], \,  W^{1, 2}(\Omega)\times W^{1, 2}(\Omega)\right), \quad \Gamma(0)=(v_{1, \lambda}, v_{2, \lambda}), \quad
\Gamma(1)=(\tilde{v}_1, \tilde{v}_2)\Big\}
\eerr
and define
 \berr
  \theta_0=\inf\limits_{\Gamma\in\mathcal{P}}\sup\limits_{t\in[0, 1]}I(\Gamma(t)).\label{c78}
 \eerr

Then we  obtain
 \be
 \theta_0> I(\tilde{v}_1, \tilde{v}_2).
  \ee

At last, noting  Lemma \ref{lem8},  \eqref{c75}--\eqref{c77}, we  can use  the mountain-pass theorem of  of Ambrosetti-Rabinowitz  \cite{amra} to conclude that
 $\theta_0$ is  also  a critical  value of the functional $I$, which gives another critical point of $I$.  In view of  \eqref{c78}, we obtain  a second solution of  \eqref{c1}--\eqref{c2},  which is different from $(v_{1, \lambda}, v_{2, \lambda})$.  Then the proof of Theorem \ref{th2} is complete.

\subsection{Quantized fluxes over $\Omega$}
 In this short subsection we calculate the quantized fluxes stated in Theorem \ref{th3} for the doubly periodic domain case.
 In fact, as in the planar case  we obtain from   \eqref{F1}--\eqref{F2}  that
\berr
F_{12}^0 &=& -\frac{1}{\sqrt{2N}} \Delta\left([N-1]\ln|\phi|^2 +|\ln\phi_N|^2\right),\\
F_{12}^{N^2-1} &=& -\sqrt{\frac{N-1}{2N}} \Delta\left(\ln|\phi|^2 - \ln\phi_N|^2\right),
 \eerr
which gives
  \berr
   {\cal F}^{U(1)}  &=&-\frac{1}{\sqrt{2N}} \ito\Delta \left([N-1]\ln|\phi|^2+|\ln\phi_N|^2\right)\ud x,\\
     {\cal F}^{SU(N)}  &=&-\sqrt{\frac{N-1}{2N}} \ito\Delta \left(\ln|\phi|^2-|\ln\phi_N|^2\right)\ud x.
  \eerr
 Then, using the equations  \eqref{a1}--\eqref{a2} and a direct integration, we get the desired quantized fluxes \eqref{d6}--\eqref{d6a}.

~

\noindent\underline{Acknowledgments}: S.C and X.H were   supported by the Natural Science Foundation of China under grant 11201118.  F.S and G.L were supported by  CONICET  , ANPCYT , CICBA, UBA and UNLP, Argentina. G.L wishes to thank LPTHE-Paris VI for hospitality and CNRS for partial financial support during part of this work.

\end{document}